\documentclass[11pt]{amsart}
\usepackage{amsmath,amssymb,amsthm}
\usepackage[margin=28mm]{geometry}
\usepackage{booktabs}
\usepackage{hyperref}

\theoremstyle{plain}
\newtheorem{theorem}{Theorem}[section]
\newtheorem{lemma}[theorem]{Lemma}
\newtheorem{proposition}[theorem]{Proposition}
\newtheorem{corollary}[theorem]{Corollary}
\newtheorem{conjecture}[theorem]{Conjecture}
\theoremstyle{definition}
\newtheorem{definition}[theorem]{Definition}

\theoremstyle{remark}
\newtheorem{remark}[theorem]{Remark}
\numberwithin{equation}{section}

\newcommand{\R}{\mathbb R}
\newcommand{\C}{\mathbb C}
\DeclareMathOperator{\Tr}{Tr}
\DeclareMathOperator{\tr}{tr}
\DeclareMathOperator{\diag}{diag}
\DeclareMathOperator{\Var}{Var}
\DeclareMathOperator{\adj}{adj}
\DeclareMathOperator{\ad}{ad}
\newcommand{\norm}[1]{\lVert #1\rVert}
\newcommand{\abs}[1]{\lvert #1\rvert}
\newcommand{\ip}[2]{\langle #1,#2\rangle}
\newcommand{\T}{\mathsf T}

\title[The weighted B\"ottcher--Wenzel inequality]{The weighted B\"ottcher--Wenzel inequality: cases of
equality and a Ky Fan norm form}
\author{Tomoharu Suzuki}
\email{quantale@gmail.com}
\urladdr{https://orcid.org/0009-0006-9163-2944}
\date{20 November 2026; revised 28 November 2026}

\begin{document}

\begin{abstract}
We prove the sharp bound $\norm{XZY-YZX}_2\le\norm Z_{(2),2}\norm X_2\norm Y_2$ for all complex matrices
of compatible sizes, where $\norm Z_{(2),2}$ is the Ky Fan $(2),2$ norm (the square root of the sum of the
squares of the two largest singular values); the constant is attained unless $X$ and $Y$ are $1\times1$
matrices.  For a positive definite matrix $\omega$ with two smallest eigenvalues $\lambda_1\le\lambda_2$ and
the weighted Frobenius norm $\norm X_\omega^2=\Tr(X^*X\omega)$, the case $Z=\omega^{-1/2}$ is the inequality
$\norm{[A,B]}_\omega^2\le\frac{\lambda_1+\lambda_2}{\lambda_1\lambda_2}\norm A_\omega^2\norm B_\omega^2$
for all real or complex matrices $A,B$ of every size, conjectured by Mayumi, Kimura, Ohno and
Chru\'sci\'nski, and we determine its cases of equality.  An earlier
preprint of Shikunov, of which the author learned after this work was completed, contains a proof of the
inequality along similar lines with a different final step.  The best constant is the square of a unitarily invariant
norm of $Z$, which is convex, so it suffices to treat weights whose $n-1$ largest
eigenvalues coincide.  For these, an inequality for one matrix and a pair of unit vectors
reduces the conjecture to the positive semidefiniteness of an explicit $4\times4$ Hermitian matrix
$F_\C(B)$, and a determinant identity, proved by hand, writes $\det F_\C$ as a sum of products of
visibly nonnegative factors; for $\lambda_1<\lambda_2$ this decomposition also yields
the cases of equality.
As a second theorem we prove, for every weight, a strengthening of the
conjecture that is equivalent to the inequality for one matrix and all pairs of unit vectors; its best
constant is again a convex function of $\omega^{-1/2}$, and the same argument applies.  We also read the inequality as a sharp bound for inner derivations on the GNS
space of a faithful state.
\end{abstract}

\maketitle

\section{Introduction}\label{sec:intro}

\subsection{The conjecture}
For $X\in M_n=M_n(\C)$ let $\norm X_2^2=\Tr(X^*X)$ be the Frobenius norm.  B\"ottcher and Wenzel
conjectured \cite{BW2005}, for real matrices, that
\begin{equation}\label{eq:BW}
\norm{XY-YX}_2^2\le2\norm X_2^2\norm Y_2^2\qquad(X,Y\in M_n).
\end{equation}
This was proved by Vong and Jin \cite{VongJin} and, independently, by Lu \cite{Lu2011}; B\"ottcher and Wenzel
\cite{BW2008} gave a different proof, which covers complex matrices.  See \cite{Lu2012} and Audenaert
\cite{Audenaert2010} for further proofs and refinements and \cite{LuWenzel} for a survey.

Let $\omega\in M_n$ be positive definite with eigenvalues $\lambda_1\le\lambda_2\le\dots\le\lambda_n$,
and put $\norm X_\omega^2=\Tr(X^*X\omega)$.  Mayumi, Kimura, Ohno and Chru\'sci\'nski \cite{MKOC24}
considered the possible ways of measuring $[A,B]$, $A$ and $B$ by $\norm\cdot_\omega$ or
$\norm\cdot_2$ in an inequality of the form \eqref{eq:BW}, determined the best constants in two of
them, and conjectured the best constants in the others.  The central one is the following.

\begin{conjecture}[{\cite[Conjecture~1, (14)]{MKOC24}}]\label{conj:MKOC}
For all $n\ge2$, all positive definite $\omega\in M_n$ and all $A,B\in M_n$,
\begin{equation}\label{eq:14}
\norm{[A,B]}_\omega^2\le\frac{\lambda_1+\lambda_2}{\lambda_1\lambda_2}\,\norm A_\omega^2\,\norm B_\omega^2 .
\end{equation}
\end{conjecture}

For $\omega=I$ this is \eqref{eq:BW}.  The constant cannot be decreased: if $e_1,e_2$ are unit
eigenvectors of $\omega$ for $\lambda_1,\lambda_2$, then $A=e_1e_2^*$ and $B=A^*$ give
$\norm{[A,B]}_\omega^2=\lambda_1+\lambda_2$, $\norm A_\omega^2=\lambda_2$ and $\norm B_\omega^2=\lambda_1$
\cite[\S2]{MKOC24}.  The same conjecture of \cite{MKOC24} has two further parts:
\begin{gather}
\norm{[A,B]}_\omega^2\le\Bigl(1+\frac{\lambda_n}{\lambda_1}\Bigr)\norm A_\omega^2\,\norm B_2^2,
\label{eq:15}\\
\norm{[A,B]}_2^2\le\frac{\lambda_1+\lambda_2}{\lambda_1^2\lambda_2}\,\norm A_\omega^2\,\norm B_\omega^2
\label{eq:16}
\end{gather}
(\cite[(15), (16)]{MKOC24}).  Inequality \eqref{eq:16} follows from \eqref{eq:14} since
$\norm C_\omega^2\ge\lambda_1\norm C_2^2$ \cite[\S2]{MKOC24}; inequality \eqref{eq:15} is a different
inequality, with the unweighted norm of $B$ on the right.  The motivation in \cite{MKOC24} is an
uncertainty relation for a state-dependent norm of a commutator \cite[\S3.1, (24)]{MKOC24}; see also
\cite{MKOC24PRA,KMO25} and Section~\ref{sec:interpretation}.

\subsection{Known results}
Table~\ref{tab:known} lists the special cases of \eqref{eq:14} that we know of.

\begin{table}[ht]
\caption{Special cases of \eqref{eq:14} in the literature (check of November 2026); for the general case see
\cite{Shikunov2026} and the text below.}
\label{tab:known}
\begin{tabular}{@{}p{0.52\textwidth}p{0.40\textwidth}@{}}
\toprule
Case & Reference\\
\midrule
$n=2$ & \cite[\S2.1]{MKOC24}\\
$A$ normal and commuting with $\omega$ & \cite[Appendix~A]{MKOC24}\\
one of $A,B$ normal, all $n$ & \cite[Supplementary Information, Theorem~1]{KMO25}\\
$A$ and $B$ both of rank one & \cite[Theorem~2.4]{FangCheng}\\
one of $A,B$ of rank one, all $n$ & \cite[Corollary~4.2]{Nicolas2026} (preprint)\\
$\lambda_1=\lambda_2$ & follows from the theorem of Shi and Zhang \cite{ShiZhang}, which proves a
conjecture of Nobori \cite{Nobori2025} (see Section~\ref{subsec:intro-main})\\
\bottomrule
\end{tabular}
\end{table}

Since Hermitian matrices are normal,
the case of Hermitian $A,B$ is contained in the third row of Table~\ref{tab:known} \cite{KMO25}.  Nicolas
\cite{Nicolas2026} also treats a factor supported on a plane invariant under $\omega$ (such a factor has rank at most
two), which reduces to the case $n=2$ \cite[Proposition~5.1]{Nicolas2026}, and for rank two he reduces \eqref{eq:14}
to the positivity of an operator pencil \cite[Propositions~6.4, 6.8]{Nicolas2026}, which he leaves open
\cite[Section~7, Questions~1--2]{Nicolas2026}.
Inequality \eqref{eq:15}
has been proved for all $n$, with a Lean formalization, in a public repository \cite{Omega14172}
(not peer reviewed; there the case \eqref{eq:14} is stated to lie outside the result), and by
Liu \cite[Theorem~1.1]{Liu2026}, who also determines the cases of equality in
\eqref{eq:15} \cite[Theorem~1.2]{Liu2026}.  Our proof of \eqref{eq:14} does not use \eqref{eq:15}.
In a preprint dated 26 September 2026 \cite{Shikunov2026} (not peer reviewed), of which we learned after the
arXiv versions v1 and v2 of this paper, Shikunov gives proofs of \eqref{eq:14}, stated in the equivalent form
with the term $\abs{\Tr(A^*B\omega)}^2$ (see the display following Theorem~\ref{cor:14}), and of \eqref{eq:15},
for all $n$ \cite[Theorem~1.1]{Shikunov2026}.  His proof of \eqref{eq:14} and ours, found independently, have
the same structure: an identity of Audenaert type and a reduction to pure states lead to an inequality for
one matrix and two unit vectors \cite[Lemmas~2.1 and~2.2, Theorem~3.1]{Shikunov2026}, which for weights whose $n-1$
largest eigenvalues coincide is our (2V); this inequality is reduced to dimension three
\cite[Section~3.1]{Shikunov2026}; and convexity and unitary invariance reduce every weight to such weights
\cite[Section~4]{Shikunov2026}.  The proofs differ in the three-dimensional case.  There Shikunov factorizes
the determinant of an explicit $9\times9$ Hermitian matrix, an identity in five real indeterminates whose
expansion is not displayed, and concludes by a continuity argument in the weight
\cite[Sections~3.2--3.4 and Appendix~A]{Shikunov2026}; we use a $4\times4$ matrix and a determinant identity
proved by hand (Theorem~\ref{thm:m2-D}).  The following are not in \cite{Shikunov2026}: the cases of equality
in \eqref{eq:14} (Theorem~\ref{thm:intro-eq}); the two-vector inequality for every weight and for arbitrary
$Z$, with the best constant (Theorem~\ref{thm:intro-2V}, Corollary~\ref{cor:m2-2V-all}); and the statement
for rectangular matrices and arbitrary $Z$ with the Ky Fan $(2),2$ norm (Theorem~\ref{thm:intro-norm}; for
positive definite $Z$ see \cite[(4.3)]{Shikunov2026}).
Apart from \cite{Shikunov2026}, we know of no case of \eqref{eq:14} without a
normality or rank assumption and with $\lambda_1<\lambda_2$, $n\ge3$, that has been settled in the
literature.  The authors of \cite{KMO25} themselves state that the weighted inequality \cite[(8)]{KMO25} was
conjectured in \cite{MKOC24} for arbitrary matrices and that they prove it when one of the two matrices is normal
(see the sentences following (22) of the journal version).

\subsection{Main results}\label{subsec:intro-main}
We prove Conjecture~\ref{conj:MKOC} in the following form, which involves no state.  Here
$s_1(Z)\ge s_2(Z)\ge\dots$ are the singular values of a matrix $Z$ (with $s_2(Z)=0$ if $Z$ has only one row
or one column), and $\norm Z_{(2),2}=(s_1(Z)^2+s_2(Z)^2)^{1/2}$ is the Ky Fan $(2),2$ norm.  We write $E_{ij}$ for
the matrix units and $X\succeq0$ ($X\succ0$) if the Hermitian matrix $X$ is positive semidefinite (definite).

\begin{theorem}[norm form]\label{thm:intro-norm}
Let $p,q\ge1$.  For all $X,Y\in M_{p\times q}(\C)$ and $Z\in M_{q\times p}(\C)$,
\begin{equation}\label{eq:intro-norm}
\norm{XZY-YZX}_2\le\norm Z_{(2),2}\,\norm X_2\,\norm Y_2 .
\end{equation}
If $(p,q)\ne(1,1)$, then for every $Z$ there are nonzero $X,Y$ with equality.  In particular, for $n\ge2$ and every
$Z\in M_n(\C)$,
\[
\mathfrak n(Z):=\sup_{X,Y\ne0}\frac{\norm{XZY-YZX}_2}{\norm X_2\norm Y_2}=\norm Z_{(2),2}\qquad(X,Y\in M_n(\C)).
\]
\end{theorem}

For $Z=I_n$ inequality \eqref{eq:intro-norm} is \eqref{eq:BW}, since $\norm{I_n}_{(2),2}^2=2$.  For a positive
definite $\omega\in M_n$ and $Z=\omega^{-1/2}$ one has $\norm Z_{(2),2}^2=\lambda_1^{-1}+\lambda_2^{-1}$, and with
$X=A\omega^{1/2}$, $Y=B\omega^{1/2}$ (so that $XZY-YZX=[A,B]\,\omega^{1/2}$, $\norm X_2=\norm A_\omega$ and
$\norm Y_2=\norm B_\omega$) inequality \eqref{eq:intro-norm} becomes \eqref{eq:14}.  Thus
Theorem~\ref{thm:intro-norm} contains the conjecture:

\begin{theorem}\label{cor:14}
Inequality \eqref{eq:14} holds for every $n\ge2$ \textup{(i)} for all real $A,B$ and every real
symmetric positive definite $\omega$, and \textup{(ii)} for all complex $A,B$ and every Hermitian
positive definite $\omega$.  The constant is best possible.
\end{theorem}

Part (i) is the special case of real matrices of part (ii); the phases of the entries of $B$ enter the proof only
through a squared modulus (Lemma~\ref{lem:m2-cross}).  The proof is short and uses no computer
(Section~\ref{subsec:intro-sketch}).

For $\lambda_1=\lambda_2$ the constant in \eqref{eq:intro-norm} with $Z=\omega^{-1/2}$ is $\sqrt2\,\norm Z_{\mathrm{op}}$
($\norm\cdot_{\mathrm{op}}$ the operator norm), and this case is contained in the theorem of Shi and Zhang
\cite{ShiZhang} (Table~\ref{tab:known}), where $\norm X_2$ is even replaced by the Ky Fan $(2),2$ norm of $X$.  It
also follows from \eqref{eq:BW} alone: $(XUY-YUX)U=[XU,YU]$ for unitary $U$, every $Z$ with $\norm Z_{\mathrm{op}}\le1$
is the mean of two unitaries (a finite-dimensional form of the theorem of Russo and Dye \cite{RussoDye}; if
$Z=U_0\diag(\cos\theta_k)V_0$, take $U_0\diag(e^{\pm i\theta_k})V_0$), and $Z\mapsto XZY-YZX$ is linear, so
$\norm{XZY-YZX}_2\le\sqrt2\norm Z_{\mathrm{op}}\norm X_2\norm Y_2$ for every $Z\in M_n$.

Since $XZY-YZX$ does not change when $X$ is replaced by $X-tY$ ($t\in\C$), inequality \eqref{eq:intro-norm} is
equivalent to
\[
\norm{XZY-YZX}_2^2\le\norm Z_{(2),2}^2\bigl(\norm X_2^2\norm Y_2^2-\abs{\Tr(X^*Y)}^2\bigr).
\]
With $\norm Z_2$ in place of $\norm Z_{(2),2}$ this is the inequality of L\'aszl\'o \cite[Theorem~3.1]{Laszlo2022}
(real square matrices) and of Nobori \cite[Proposition~3.1]{Nobori2025} (complex rectangular matrices), and with
$\sqrt2\,\norm Z_{\mathrm{op}}$ in place of $\norm Z_{(2),2}$ it is \cite[(9)]{Nobori2025}, which Nobori derived from
his conjecture proved by Shi and Zhang \cite{ShiZhang}.  Since
$\norm Z_{(2),2}\le\min\{\norm Z_2,\sqrt2\,\norm Z_{\mathrm{op}}\}$, Theorem~\ref{thm:intro-norm} sharpens both, and the
factor $\norm Z_{(2),2}$ is best possible for every $Z$.

We also determine the cases of equality in \eqref{eq:14}.

\begin{theorem}[equality]\label{thm:intro-eq}
Let $\omega\in M_n$ be positive definite, $\varphi=\Tr(\omega\,\cdot\,)$, and $A,B\ne0$.  Equality holds in
\eqref{eq:14} if and only if there is a two-dimensional subspace $P$, invariant under $\omega$, on which
$\omega$ has the eigenvalues $\lambda_1,\lambda_2$ (counted with multiplicity), such that the orthogonal projection $\Pi$ onto $P$ satisfies
$A=\Pi A\Pi$, $B=\Pi B\Pi$ and $\varphi(A)=\varphi(B)=\varphi(A^*B)=0$; that is, $\Pi$, $A$ and $B$ are mutually
orthogonal for the inner product $\varphi(X^*Y)$.  For real $A,B$ and real symmetric $\omega$, every such $P$ is real.
\end{theorem}

In an eigenbasis, $A=e_1e_2^*$, $B=A^*$ (the example above) and $A=\lambda_2E_{11}-\lambda_1E_{22}$, $B=E_{12}$
are such pairs.  For $\omega=I$ the theorem is the description of the pairs with equality in \eqref{eq:BW}
by Cheng, Vong and Wenzel \cite[Theorem~3.1]{CVW2010} (see \cite[Theorem~2.7]{LuWenzel}), which our proof uses when
$\lambda_1=\lambda_2$; for $\lambda_1<\lambda_2$ it uses only the tools of the proof of Theorem~\ref{cor:14}
(Section~\ref{subsec:m2-equality}).  The same proof gives the cases of equality in \eqref{eq:intro-norm} for every $Z$
with $s_2(Z)>0$ (Remark~\ref{rem:eq-Z}).  For inequality \eqref{eq:15} with $\lambda_1<\lambda_n$, the analogous
description, with a rank-one nilpotent second factor supported on a plane spanned by eigenvectors for $\lambda_1$
and $\lambda_n$, is due to Liu \cite[Theorem~1.2]{Liu2026}.  Since \eqref{eq:16} follows from \eqref{eq:14}, we also obtain:

\begin{corollary}\label{cor:16}
Inequality \eqref{eq:16} holds for every $n\ge2$, all complex $A,B$ and every Hermitian positive
definite $\omega$.
\end{corollary}
For Hermitian $A,B$ a sharper constant, $2/(\lambda_1\lambda_2)$, was obtained in
\cite[Corollary~D(iii)]{SuzukiP12}; see Remark~\ref{rem:int-P12}.  If $\lambda_1<\lambda_2$ and $A,B$ both have rank
one, the constant in \eqref{eq:16} is not best possible; the best constant in this case, which depends only on
$\lambda_1$ and $\lambda_2$, and the pairs attaining it were determined by Pang \cite[Theorems~1.1 and~5.1,
Corollary~5.2]{Pang2026}.

Theorem~\ref{cor:14} answers Questions~1 and~2 of \cite[Section~7]{Nicolas2026} in the affirmative
(as does \cite[Theorem~1.1]{Shikunov2026}).
Applying \eqref{eq:red-TD} to $\sqrt t\,R_1+t^{-1/2}R_2$ ($t>0$), where $R_1,R_2$ are the two
orthogonal rank-one singular directions of Frobenius norm one, gives the positivity of his operator
pencil for every $t$ \cite[Proposition~6.4]{Nicolas2026}, and by \cite[Proposition~6.8]{Nicolas2026}
(when $R_1$ attains equality in the rank-one bound, as assumed there) this implies the positive
semidefiniteness of his tangency operator.
Theorem~\ref{cor:14} also settles the rank-two case and the case $n=3$ \cite[Section~7, Question~7]{Nicolas2026}.
Our route is different: Theorem~\ref{thm:m2-D}, on which our proof rests, does not use the pencil of
\cite{Nicolas2026}.

\subsection{A second theorem: the two-vector inequality}\label{subsec:intro-second}
The proof of Theorem~\ref{cor:14} passes through an inequality for one matrix and a pair of unit
vectors, the inequality (2V) of Definition~\ref{def:red-2V}, which is stronger than \eqref{eq:14}; the
convexity argument of the proof applies to it verbatim.  Our second theorem is
(2V) for every weight, in the following equivalent form.

\begin{theorem}[two-vector inequality]\label{thm:intro-2V}
Let $\omega\in M_n$ be positive definite, $D=\omega^{-1/2}$ and $c=\lambda_1^{-1}+\lambda_2^{-1}$.  For all
$X,Y\in M_n$ with $Y\ne0$, with $\mathsf W:=DX^*Y+YX^*D$,
\begin{equation}\label{eq:intro-2V}
\norm{XDY-YDX}_2^2+\norm{\mathsf W}_2^2-\frac{\abs{\Tr(Y^*\mathsf W)}^2}{\norm Y_2^2}\le c\,\norm X_2^2\norm Y_2^2 .
\end{equation}
\end{theorem}

By the Cauchy--Schwarz inequality the last two terms on the left add up to a nonnegative number, so
with $X=A\omega^{1/2}$, $Y=B\omega^{1/2}$ Theorem~\ref{thm:intro-2V} implies \eqref{eq:14}, and for each pair it
is stronger.  For $Y=\psi_1\psi_2^*$ with unit vectors $\psi_1,\psi_2$, \eqref{eq:intro-2V} is (2V) for the pair
$(\psi_1,\psi_2)$ (Lemma~\ref{lem:red-pure}(a)), and (2V) for all pairs implies \eqref{eq:intro-2V}
(Lemma~\ref{lem:red-pure}).  With $Z^*X^*Y+YX^*Z^*$ in place of $\mathsf W$ and $\norm Z_{(2),2}^2$ in place of $c$,
\eqref{eq:intro-2V} holds for every $Z\in M_n(\C)$ in place of $D$, with the best constant
(Corollary~\ref{cor:m2-2V-all}); the adjoint is needed (Remark~\ref{rem:red-two-level-2V}).  On the GNS space of
$\varphi$, \eqref{eq:intro-2V} adds to the commutator a modular-twisted anticommutator, and for $\omega=I$ it is the
variance bound of Audenaert \cite{Audenaert2010} (Proposition~\ref{prop:int-2V-gns} and Remark~\ref{rem:int-2V-gns}).

Both theorems rest on the positive semidefiniteness of an explicit $4\times4$ matrix $F_\C(B)$ for $2\times2$
matrices $B$ in a convex set $K^\C_2$ (Definition~\ref{def:F}, Theorem~\ref{thm:m2-D}).  The analogous
inequality for $m\times m$ matrices $B$, $m\ge3$, is a matrix inequality that we treat in a separate paper (in preparation); the
proofs in this paper need only $m=2$.

\subsection{Sketch of the proof}\label{subsec:intro-sketch}
With $X=A\omega^{1/2}$, $Y=B\omega^{1/2}$ and $D=\omega^{-1/2}=\diag(d_1\ge d_2\ge\dots\ge d_n)$ in an eigenbasis
of $\omega$, \eqref{eq:14} becomes
\[
\norm{T_D(X,Y)}_2^2\le c\,\norm X_2^2\norm Y_2^2,\qquad T_D(X,Y):=XDY-YDX,\qquad
c:=d_1^2+d_2^2=\frac1{\lambda_1}+\frac1{\lambda_2}
\]
(Lemma~\ref{lem:red-normal}); the same substitution is used by Nicolas \cite{Nicolas2026}.  The proof has four
steps.

\emph{Linearisation in one factor.}  We fix $X$.  A weighted form of an identity of Audenaert
\cite{Audenaert2010} (Lemma~\ref{lem:red-Aud}) writes $\norm{T_D(X,Y)}_2^2$ as a linear functional of the pair of
states $(YY^*,Y^*Y)/\norm Y_2^2$, times $\norm Y_2^2$, minus the squared norm of a matrix depending linearly on $Y$.  The
Cauchy--Schwarz inequality and a reduction to pure pairs (by the Toeplitz--Hausdorff theorem;
Lemma~\ref{lem:red-pure}) reduce both theorems to the inequality (2V), for one matrix $X$ and two unit vectors
$\psi_1,\psi_2$.  That $T_D$ is evaluated at pairs $(X,Y)$ is essential here (Remark~\ref{rem:intro-wedge}).

\emph{Convexity.}  The first step holds for every weight; the second step shows that it is enough to treat weights
with $\lambda_2=\dots=\lambda_n$.  The best constant is $\mathfrak n(D)^2$, where $\mathfrak n$ (Theorem~\ref{thm:intro-norm}) is a
seminorm on $M_n$ with $\mathfrak n(UZV)=\mathfrak n(Z)$ for unitaries $U,V$, and $c=\norm D_{(2),2}^2$.  Writing a
diagonal matrix explicitly as a convex combination of matrices with singular values $(a,b,\dots,b)$ (which, normalised,
are the extreme points of the unit ball of the Ky Fan $(2),2$ norm on diagonal matrices) reduces
Theorem~\ref{thm:intro-norm}, and with it \eqref{eq:14}, to
weights with $\lambda_2=\dots=\lambda_n$ (Proposition~\ref{prop:red-two-level}).  The same applies to the best constant
in \eqref{eq:intro-2V}, with $D^*$ in place of $D$ in $\mathsf W$ for non-Hermitian arguments.

\emph{Localization.}  For $D=bI+(a-b)ee^*$ with a unit vector $e$, every subspace containing $e$ is invariant under $D$,
and $D=b$ on its orthogonal complement.  For a three-dimensional such subspace $V$ containing $\psi_1,\psi_2$, the
functional of (2V) splits into its value at the compression of $X$ to $V$ and two terms bounded by $c$ times the
squared norms of the off-diagonal blocks of $X$ (Lemma~\ref{lem:red-local}); so (2V) follows from the case $n=3$
(Corollary~\ref{cor:red-local}).

\emph{The $4\times4$ matrix.}  For $n=3$, in bases adapted to $\psi_1,\psi_2$ (Definition~\ref{def:red-adapted}), the
corner of $D$ is a matrix $B\in M_2(\C)$, and $R,S,H$ are the normalised squared moduli of the three other nonzero
entries.  Completing squares (Theorem~\ref{thm:red-square}) leaves the Hermitian form of the $4\times4$ matrix
$F_\C(B)$ (Definition~\ref{def:F}), and the Ky Fan maximum principle shows $B\in K^\C_2$
(Lemma~\ref{lem:red-KyFan}).  The positive semidefiniteness of $F_\C(B)$ (Theorem~\ref{thm:m2-D}) is proved by hand:
an identity writes $\det F$, up to a positive factor, as a product of two nonnegative quantities plus $HRS$ times
a sum of products of nonnegative factors, and since $F$ is a positive definite matrix minus a matrix of rank one,
$\det F>0$ gives $F\succ0$ in the interior of $K^\C_2$.

The assembly is Section~\ref{subsec:m2-proof14}, and the cases of equality are determined in
Section~\ref{subsec:m2-equality} with the same tools.  Section~\ref{sec:interpretation} reads \eqref{eq:14} in terms
of the GNS representation of a state and of uncertainty relations.

\begin{remark}\label{rem:intro-wedge}
The map $T_D$ factors through a linear map $\tilde T_D$ on ${\wedge}^2(M_n)$, normed by
$\norm{X\wedge Y}^2=\norm X_2^2\norm Y_2^2-\abs{\Tr(X^*Y)}^2$, but $c$ is not a bound for its norm: for $\omega=I$
(so $c=2$) and $\xi=E_{12}\wedge E_{21}+E_{13}\wedge E_{31}$ one has $\norm\xi^2=2$ and
$\tilde T_I(\xi)=2E_{11}-E_{22}-E_{33}$, with squared norm $6=3\norm\xi^2$.  On decomposable elements it is a bound:
$\norm{T_D(X,Y)}_2^2\le c\,\norm{X\wedge Y}^2$, by Theorem~\ref{thm:intro-norm} with $Z=D$ applied to $X-tY$ and $Y$,
where $t=\Tr(Y^*X)/\norm Y_2^2$ ($Y\ne0$).
\end{remark}

\begin{remark}[an inductive approach]\label{rem:H}
A natural approach to \eqref{eq:14} is induction on $n$: for $D=\diag(d,1,\dots,1)$, $d>1$, one would like every
critical point $(X,Y)$ of $\norm{XDY-YDX}_2^2/\norm{X\wedge Y}^2$ (Remark~\ref{rem:intro-wedge}) whose support (the
sum of the ranges of $X,Y,X^*,Y^*$) is all of $\C^n$ to admit a direction $X\mapsto X+t(AX+XA)$, $Y\mapsto Y+t(AY+YA)$, with $A$
Hermitian and commuting with $D$, along which the ratio increases to second order (statement (H); it would let
the induction proceed).  Statement (H)
fails for real matrices with $n=4$: for $1<d<2$ there is a real critical point, supported on all of $\R^4$ and with
no common invariant decomposition, along which no real direction of this form increases the ratio to second order;
its value $2/(3-d)$ is below $d^2+1$.  The point will be given in the separate paper (in preparation)
mentioned in Section~\ref{subsec:intro-second}.  We do not know whether the complex form of (H) holds.
\end{remark}

\subsection{Role of the computer}\label{subsec:intro-computer}
This paper uses no computer-assisted proof: all identities, including the coefficient table of the case $m=2$,
are proved by hand.  Symbolic and numerical checks, with code written using Claude (Anthropic) and checked by the
author, were used only as an independent confirmation.

\section{From the conjecture to a matrix inequality}\label{sec:reduction}

In this section we reduce \eqref{eq:14} and Theorem~\ref{thm:intro-2V} to an inequality (2V) for one
matrix and two unit vectors (Section~\ref{subsec:red-2V}), the weights we need to $\lambda_2=\dots=\lambda_n$ and then
to $n=3$ (Section~\ref{subsec:red-two-level}), and (2V) for $n=3$ to the positive semidefiniteness of the $4\times4$
matrices $F_\C(B)$ of Definition~\ref{def:F} with $m=2$ (Sections~\ref{subsec:red-adapted}
to~\ref{subsec:red-assembly}):
\[
\text{\eqref{eq:14}}\ \Longleftarrow\ \text{(2V)}\ \Longleftarrow\ \text{(2V) for $n=3$}\ \Longleftarrow\
F_\C(B)\succeq0\ \text{on }K^\C_2 ,
\]
where the second implication holds for weights with $\lambda_2=\dots=\lambda_n$, to which
Proposition~\ref{prop:red-two-level} reduces.
We write every step for complex matrices and Hermitian $\omega$;
real $A,B$ and real symmetric $\omega$ are a special case.

\emph{Conventions.}  $\langle a,b\rangle=a^*b$ is conjugate linear in the first variable,
$\langle X,Y\rangle=\Tr(X^*Y)$ for matrices, and $\norm X=\norm X_2$ is the Frobenius norm (the Euclidean
norm on vectors).  $\bar X$ is the
entrywise conjugate, $X^\T$ the transpose and $X^*=\bar X^\T$; $I_k$ is the $k\times k$ identity matrix, and
$(x;y)$ is the column vector obtained by stacking $x$ above $y$.  A \emph{state} is a positive
semidefinite Hermitian matrix of trace one; $\mathfrak S_n$ denotes the set of $n\times n$ states, a
compact convex subset of the real vector space of Hermitian matrices.  For a unit vector $\psi$ we
write $P^\perp_\psi=I-\psi\psi^*$.

\subsection{Normalisation}

Unless stated otherwise, in this section $D=\diag(d_1\ge d_2\ge\cdots\ge d_n>0)$, $c=d_1^2+d_2^2$ and
$T_D(X,Y)=XDY-YDX$ for $X,Y\in M_n(\C)$.

\begin{lemma}[normalisation]\label{lem:red-normal}
Let $\omega$ be a positive definite Hermitian $n\times n$ matrix, $Z=\omega^{-1/2}$, and for $A,B\in M_n(\C)$ put
$X=A\omega^{1/2}$, $Y=B\omega^{1/2}$ \textup{(}a bijection\textup{)}.  Then $\norm X=\norm A_\omega$,
$\norm Y=\norm B_\omega$ and $XZY-YZX=[A,B]\,\omega^{1/2}$, so $\norm{XZY-YZX}=\norm{[A,B]}_\omega$.  In particular, if
$\omega=\diag(\lambda_1(\omega)\le\lambda_2(\omega)\le\cdots)$ is diagonal and $D=\omega^{-1/2}$, then
\eqref{eq:14} holds for $\omega$ and all complex $A,B$ if and only if
\begin{equation}\label{eq:red-TD}
\norm{T_D(X,Y)}^2\le c\,\norm X^2\norm Y^2\qquad\text{for all }X,Y\in M_n(\C).
\end{equation}
\end{lemma}

\begin{proof}
$\norm{C\omega^{1/2}}^2=\Tr(\omega^{1/2}C^*C\omega^{1/2})=\Tr(C^*C\omega)=\norm C_\omega^2$ for every $C$, and
$XZY=AB\omega^{1/2}$, $YZX=BA\omega^{1/2}$.  For diagonal $\omega$, $T_D(X,Y)=XZY-YZX$ and
$c=\frac1{\lambda_1}+\frac1{\lambda_2}=\frac{\lambda_1+\lambda_2}{\lambda_1\lambda_2}$.
\end{proof}

\subsection{The two-vector inequality implies the conjecture}\label{subsec:red-2V}

The definitions of this subsection and Lemmas~\ref{lem:red-Aud}, \ref{lem:red-pure} and~\ref{lem:red-decomp} use
only that $D$ is Hermitian; in Section~\ref{subsec:red-two-level} they are applied to positive definite matrices $D$
that need not be diagonal, with $c$ the sum of the squares of the two largest eigenvalues of $D$.

Fix $X\in M_n(\C)$ and put
\[
\mathsf M_1=DX^*XD+XD^2X^*,\qquad \mathsf M_2=DXX^*D+X^*D^2X
\]
(both Hermitian).
For a pair of states $\varpi=(\varpi_1,\varpi_2)$ put
\[
\chi(\varpi)=\Tr(DX^*\varpi_1)+\Tr(X^*D\varpi_2)\in\C,\qquad
\Upsilon(\varpi)=\Tr(\mathsf M_1\varpi_1)+\Tr(\mathsf M_2\varpi_2)-\abs{\chi(\varpi)}^2\in\R ,
\]
and for unit vectors $\psi_1,\psi_2$ write
$\Upsilon(\psi_1,\psi_2)=\Upsilon(\psi_1\psi_1^*,\psi_2\psi_2^*)$.

\begin{definition}\label{def:red-2V}
Fix $D$ and unit vectors $\psi_1,\psi_2\in\C^n$.  For $\Lambda>0$, \emph{(2V)$[\Lambda]$} is the
statement
\[
\Upsilon(\psi_1,\psi_2)\le\Lambda\norm X^2\qquad\text{for all }X\in M_n(\C).
\]
\emph{(2V)} is the statement that (2V)$[c]$ holds for all unit vectors $\psi_1,\psi_2$.
\end{definition}

\begin{lemma}[weighted Audenaert identity]\label{lem:red-Aud}
For all $X,Y\in M_n(\C)$, with $\mathsf W=DX^*Y+YX^*D$,
\[
\norm{T_D(X,Y)}^2+\norm{\mathsf W}^2=\Tr(\mathsf M_1YY^*)+\Tr(\mathsf M_2Y^*Y).
\]
\end{lemma}

\begin{proof}
Expanding,
\begin{align*}
\norm{T_D(X,Y)}^2&=\norm{XDY}^2+\norm{YDX}^2-2\operatorname{Re}\Tr(Y^*DX^*YDX),\\
\norm{\mathsf W}^2&=\norm{DX^*Y}^2+\norm{YX^*D}^2+2\operatorname{Re}\Tr(Y^*XD\,YX^*D).
\end{align*}
The cross terms cancel, because $\overline{\Tr(Y^*DX^*YDX)}=\Tr(X^*DY^*XDY)=\Tr(Y^*XDYX^*D)$ by
cyclicity.  For the square terms,
\[
\norm{XDY}^2+\norm{DX^*Y}^2=\Tr(\mathsf M_1YY^*),\qquad \norm{YDX}^2+\norm{YX^*D}^2=\Tr(\mathsf M_2Y^*Y),
\]
since by cyclicity
\begin{align*}
\norm{XDY}^2&=\Tr(DX^*XD\,YY^*),&\norm{DX^*Y}^2&=\Tr(XD^2X^*\,YY^*),\\
\norm{YDX}^2&=\Tr(DXX^*D\,Y^*Y),&\norm{YX^*D}^2&=\Tr(X^*D^2X\,Y^*Y).
\end{align*}
\end{proof}

\begin{lemma}[reduction to pure pairs]\label{lem:red-pure}
Fix $X\in M_n(\C)$.
\begin{enumerate}
\item[(a)] If $Y\ne0$, $\varpi_1=YY^*/\norm Y^2$ and $\varpi_2=Y^*Y/\norm Y^2$, then $(\varpi_1,\varpi_2)$
is a pair of states, $\langle Y,\mathsf W\rangle=\norm Y^2\chi(\varpi)$ with $\mathsf W=DX^*Y+YX^*D$, and
\[
\norm{T_D(X,Y)}^2+\norm{\mathsf W}^2-\frac{\abs{\langle Y,\mathsf W\rangle}^2}{\norm Y^2}=\norm Y^2\,\Upsilon(\varpi_1,\varpi_2);
\]
in particular $\norm{T_D(X,Y)}^2\le\norm Y^2\,\Upsilon(\varpi_1,\varpi_2)$, by the Cauchy--Schwarz inequality.
\item[(b)] $\sup_{\varpi\in\mathfrak S_n\times\mathfrak S_n}\Upsilon(\varpi)
=\max\{\Upsilon(\psi_1,\psi_2):\psi_1,\psi_2\in\C^n\text{ unit vectors}\}$.
\end{enumerate}
\end{lemma}

\begin{proof}
(a) We have
$\langle Y,\mathsf W\rangle=\Tr(Y^*DX^*Y)+\Tr(Y^*YX^*D)=\Tr(DX^*\,YY^*)+\Tr(X^*D\,Y^*Y)
=\norm Y^2\chi(\varpi)$.  By Lemma~\ref{lem:red-Aud},
$\norm{T_D(X,Y)}^2+\norm{\mathsf W}^2=\norm Y^2\bigl(\Tr(\mathsf M_1\varpi_1)+\Tr(\mathsf M_2\varpi_2)\bigr)$; subtracting
$\abs{\langle Y,\mathsf W\rangle}^2/\norm Y^2=\norm Y^2\abs{\chi(\varpi)}^2$ gives the identity, and
$\abs{\langle Y,\mathsf W\rangle}^2\le\norm Y^2\norm{\mathsf W}^2$ gives the inequality.

(b) \emph{Step 1: a maximiser.}  $\mathfrak S_n\times\mathfrak S_n$ is compact and $\Upsilon$ is
continuous, so there is a maximiser $\varpi^0=(\varpi^0_1,\varpi^0_2)$.  Put $\mu=\chi(\varpi^0)\in\C$.

\emph{Step 2: an affine majorant.}  From the scalar identity
$-\abs\chi^2=\abs\mu^2-2\operatorname{Re}(\bar\mu\chi)-\abs{\chi-\mu}^2$ and
$2\operatorname{Re}\bigl(\bar\mu\Tr(DX^*\varpi_1)\bigr)=\Tr\bigl((\bar\mu DX^*+\mu XD)\varpi_1\bigr)$
(because $\overline{\Tr(DX^*\varpi_1)}=\Tr(\varpi_1XD)$ for Hermitian $\varpi_1$), and the analogous
formula for $\varpi_2$, we obtain
\[
\Upsilon(\varpi)=L(\varpi)-\abs{\chi(\varpi)-\mu}^2,\qquad
L(\varpi)=\Tr(\mathsf H_1\varpi_1)+\Tr(\mathsf H_2\varpi_2)+\abs\mu^2,
\]
with the Hermitian matrices $\mathsf H_1=\mathsf M_1-\bar\mu DX^*-\mu XD$ and
$\mathsf H_2=\mathsf M_2-\bar\mu X^*D-\mu DX$.  The function $L$ is real affine in $\varpi$, and
$L(\varpi^0)=\Upsilon(\varpi^0)$.

\emph{Step 3: first-order condition.}  For any $\varpi$ and $t\in(0,1]$ the pair
$\varpi^t=\varpi^0+t(\varpi-\varpi^0)$ lies in $\mathfrak S_n\times\mathfrak S_n$.  Since $\chi$ and $L$ are
affine,
$\Upsilon(\varpi^t)=L(\varpi^0)+t\bigl(L(\varpi)-L(\varpi^0)\bigr)-t^2\abs{\chi(\varpi)-\mu}^2$.  From
$\Upsilon(\varpi^t)\le\Upsilon(\varpi^0)=L(\varpi^0)$ we get
$t(L(\varpi)-L(\varpi^0))\le t^2\abs{\chi(\varpi)-\mu}^2$; dividing by $t$ and letting $t\downarrow0$
gives $L(\varpi)\le L(\varpi^0)$.  Thus $\varpi^0$ maximises $L$ (no concavity of $\Upsilon$ is used).

\emph{Step 4: support.}  $L$ is a sum of a term in $\varpi_1$ and a term in $\varpi_2$, so $\varpi^0_1$
maximises $\Tr(\mathsf H_1\varpi_1)$ over $\mathfrak S_n$.  This maximum is
$\lambda_{\max}(\mathsf H_1)$, hence $\Tr\bigl((\lambda_{\max}(\mathsf H_1)I-\mathsf H_1)\varpi^0_1\bigr)=0$.
For positive semidefinite $P,Q$ one has $\Tr(PQ)=\Tr(P^{1/2}QP^{1/2})=\norm{Q^{1/2}P^{1/2}}^2$, so
$\Tr(PQ)=0$ forces $Q^{1/2}P^{1/2}=0$ and hence $QP=Q^{1/2}(Q^{1/2}P^{1/2})P^{1/2}=0$.  Applied to
$P=\lambda_{\max}(\mathsf H_1)I-\mathsf H_1$ and $Q=\varpi^0_1$ this gives
$\operatorname{range}\varpi^0_1\subseteq\mathcal T_1:=\ker(\lambda_{\max}(\mathsf H_1)I-\mathsf H_1)$.
Likewise $\operatorname{range}\varpi^0_2\subseteq\mathcal T_2:=\ker(\lambda_{\max}(\mathsf H_2)I-\mathsf H_2)$.

\emph{Step 5: choice of the pure pair.}  Let $P_1$ be the orthogonal projection onto
$\mathcal T_1$ and regard $C_1=P_1DX^*P_1$ as an operator on $\mathcal T_1$.  Writing
$\varpi^0_1=\sum_jt_j\phi_j\phi_j^*$ with $t_j\ge0$, $\sum_jt_j=1$ and unit vectors
$\phi_j\in\mathcal T_1$, we get
$\Tr(DX^*\varpi^0_1)=\sum_jt_j\langle\phi_j,C_1\phi_j\rangle$, a convex combination of points of the
numerical range $W(C_1)=\{\langle\psi,C_1\psi\rangle:\psi\in\mathcal T_1,\ \norm\psi=1\}$.  By the
Toeplitz--Hausdorff theorem \cite[\S I.2, Problem~I.6.2]{BhatiaMA} $W(C_1)$ is convex, so there is a unit vector
$\psi_1\in\mathcal T_1$ with $\langle\psi_1,DX^*\psi_1\rangle=\Tr(DX^*\varpi^0_1)$ (if
$\dim\mathcal T_1=1$, simply $\varpi^0_1=\psi_1\psi_1^*$).  Independently, choose a unit vector $\psi_2\in\mathcal T_2$ with
$\langle\psi_2,X^*D\psi_2\rangle=\Tr(X^*D\varpi^0_2)$ in the same way.

\emph{Step 6: value.}  Now $\chi(\psi_1,\psi_2)=\Tr(DX^*\varpi^0_1)+\Tr(X^*D\varpi^0_2)=\mu$, and
$\Tr(\mathsf H_k\psi_k\psi_k^*)=\lambda_{\max}(\mathsf H_k)$ since $\psi_k\in\mathcal T_k$.  Hence
$\Upsilon(\psi_1,\psi_2)=L(\psi_1,\psi_2)
=\lambda_{\max}(\mathsf H_1)+\lambda_{\max}(\mathsf H_2)+\abs\mu^2=L(\varpi^0)=\Upsilon(\varpi^0)$.  The
reverse inequality is trivial, since pure pairs are pairs of states.
\end{proof}

\begin{proposition}[(2V) implies the inequalities for $D$]\label{prop:red-2V}
If (2V) holds for $D$, then for all $X,Y\in M_n(\C)$ with $Y\ne0$, with $\mathsf W=DX^*Y+YX^*D$,
\[
\norm{T_D(X,Y)}^2+\norm{\mathsf W}^2-\frac{\abs{\langle Y,\mathsf W\rangle}^2}{\norm Y^2}\le c\,\norm X^2\norm Y^2 ;
\]
in particular \eqref{eq:red-TD} holds.
\end{proposition}

\begin{proof}
By Lemma~\ref{lem:red-pure}(a),(b) and (2V) the left side equals
\[
\norm Y^2\Upsilon(\varpi_1,\varpi_2)\le\norm Y^2\max_\psi\Upsilon(\psi_1,\psi_2)\le c\norm X^2\norm Y^2 .
\]
For $Y=0$, \eqref{eq:red-TD} is trivial.
\end{proof}

\begin{lemma}[decomposition at a pure pair]\label{lem:red-decomp}
For unit vectors $\psi_1,\psi_2$,
\[
\Upsilon(\psi_1,\psi_2)=\norm{T_D(X,\psi_1\psi_2^*)}^2+\norm{P^\perp_{\psi_1}DX^*\psi_1}^2
+\norm{P^\perp_{\psi_2}DX\psi_2}^2 .
\]
\end{lemma}

\begin{proof}
The matrix $Y=\psi_1\psi_2^*$ has $\norm Y=1$, $YY^*=\psi_1\psi_1^*$ and $Y^*Y=\psi_2\psi_2^*$.  By
Lemma~\ref{lem:red-pure}(a),
$\Upsilon(\psi_1,\psi_2)=\norm{T_D(X,Y)}^2+\norm{\mathsf W}^2-\abs{\langle Y,\mathsf W\rangle}^2
=\norm{T_D(X,Y)}^2+\norm{\mathsf W-\langle Y,\mathsf W\rangle Y}^2$.
With $\mathbf b_1=DX^*\psi_1$ and $\mathbf b_2=DX\psi_2$ we have
$\mathsf W=\mathbf b_1\psi_2^*+\psi_1\mathbf b_2^*$.  Split
$\mathbf b_1=\langle\psi_1,\mathbf b_1\rangle\psi_1+\mathbf b_1^\perp$ and
$\mathbf b_2=\langle\psi_2,\mathbf b_2\rangle\psi_2+\mathbf b_2^\perp$; then
$\mathsf W=\bigl(\langle\psi_1,\mathbf b_1\rangle+\langle\mathbf b_2,\psi_2\rangle\bigr)\psi_1\psi_2^*
+\mathbf b_1^\perp\psi_2^*+\psi_1(\mathbf b_2^\perp)^*$, and the three terms are mutually orthogonal
because $\mathbf b_1^\perp\perp\psi_1$ and $\mathbf b_2^\perp\perp\psi_2$.  Hence the first term is
$\langle Y,\mathsf W\rangle Y$ and
$\norm{\mathsf W-\langle Y,\mathsf W\rangle Y}^2=\norm{\mathbf b_1^\perp}^2+\norm{\mathbf b_2^\perp}^2$.
\end{proof}

\subsection{Weights with \texorpdfstring{$\lambda_2=\dots=\lambda_n$}{lambda2=...=lambdan}}\label{subsec:red-two-level}

In this subsection $n\ge2$ is arbitrary.
For $Z\in M_n(\C)$ let $s_1(Z)\ge s_2(Z)\ge\dots\ge s_n(Z)\ge0$ be its singular values,
$\norm Z_{(2),2}=(s_1(Z)^2+s_2(Z)^2)^{1/2}$ and $T_Z(X,Y)=XZY-YZX$.  For $X,Y\in M_n(\C)$ with $Y\ne0$ put
\[
\mathsf W_Z(X,Y)=Z^*X^*Y+YX^*Z^*,\qquad
\mathcal L_Z(X,Y)=\norm{T_Z(X,Y)}^2+\norm{P^\perp_Y\mathsf W_Z(X,Y)}^2 ,
\]
where $P^\perp_Y$ is the orthogonal projection of $M_n(\C)$ (with $\ip XY=\Tr(X^*Y)$) onto the orthogonal
complement of $Y$, so that $\norm{P^\perp_Y\mathsf W}^2=\norm{\mathsf W}^2-\abs{\ip Y{\mathsf W}}^2/\norm Y^2$, and
\[
\mathfrak n(Z)=\sup_{X,Y\ne0}\frac{\norm{T_Z(X,Y)}}{\norm X\norm Y},\qquad
\phi(Z)=\sup_{X,Y\ne0}\frac{\mathcal L_Z(X,Y)^{1/2}}{\norm X\norm Y}\qquad(X,Y\in M_n(\C)).
\]
For Hermitian $Z$, $\mathsf W_Z(X,Y)=ZX^*Y+YX^*Z$; for $Z=D$, $\mathcal L_D(X,Y)$ is the left side of the inequality
of Proposition~\ref{prop:red-2V}, and $\norm D_{(2),2}^2=c$.

\begin{lemma}[averaging over signs]\label{lem:red-sign}
Let $s\in\R^n$ and $\beta>0$ with $\abs{s_k}\le\beta$ for $k\ge2$.  For
$\epsilon=(\epsilon_2,\dots,\epsilon_n)\in\{\pm1\}^{n-1}$ put
\[
D_\epsilon=\diag(s_1,\epsilon_2\beta,\dots,\epsilon_n\beta),\qquad
w_\epsilon=\prod_{k=2}^n\frac12\Bigl(1+\frac{\epsilon_ks_k}\beta\Bigr).
\]
Then $w_\epsilon\ge0$, $\sum_\epsilon w_\epsilon=1$ and $\sum_\epsilon w_\epsilon D_\epsilon=\diag(s)$.  Moreover
$w_\epsilon>0$ if and only if $\epsilon_ks_k\ne-\beta$ for every $k\ge2$.
\end{lemma}

\begin{proof}
Each factor of $w_\epsilon$ lies in $[0,1]$ since $\abs{s_k}\le\beta$, and it vanishes exactly when
$\epsilon_ks_k=-\beta$.  Summing the product over $\epsilon$ factor by factor gives
$\sum_\epsilon w_\epsilon=\prod_{k\ge2}\bigl(\frac12(1+\frac{s_k}\beta)+\frac12(1-\frac{s_k}\beta)\bigr)=1$ and, for
$k\ge2$, $\sum_\epsilon w_\epsilon\epsilon_k\beta=\beta\bigl(\frac12(1+\frac{s_k}\beta)-\frac12(1-\frac{s_k}\beta)\bigr)=s_k$.
The first diagonal entry of every $D_\epsilon$ is $s_1$.
\end{proof}

\begin{proposition}[reduction to $\lambda_2=\dots=\lambda_n$]\label{prop:red-two-level}
Let $n\ge2$ and let $\theta$ be $\mathfrak n$ or $\phi$.  Then $\theta$ is a seminorm on $M_n(\C)$ with
$\theta(UZV)=\theta(Z)$ for all unitaries $U,V$, and $\theta(Z)\ge\mathfrak n(Z)\ge\norm Z_{(2),2}$ for all $Z$.  If
$\theta(Z)\le\norm Z_{(2),2}$ for every $Z=\diag(a,b,\dots,b)$ with $a\ge b>0$, then $\theta(Z)=\norm Z_{(2),2}$ for every
$Z\in M_n(\C)$.
\end{proposition}

\begin{proof}
\begin{enumerate}
\item \emph{Seminorm and invariance.}  Fix $X,Y\ne0$.  The map $Z\mapsto T_Z(X,Y)$ is complex linear and
$Z\mapsto P^\perp_Y\mathsf W_Z(X,Y)$ is conjugate linear, so $Z\mapsto(T_Z(X,Y),P^\perp_Y\mathsf W_Z(X,Y))$ is a real
linear map into $M_n(\C)^2$ which multiplies norms by $\abs\alpha$ when $Z$ is multiplied by $\alpha\in\C$.  Hence
$\norm{T_Z(X,Y)}$ and $\mathcal L_Z(X,Y)^{1/2}$ are seminorms in $Z$.  Since $\norm{T_Z(X,Y)}\le2\norm Z_{\mathrm{op}}\norm X\norm Y$
and $\norm{\mathsf W_Z(X,Y)}\le2\norm Z_{\mathrm{op}}\norm X\norm Y$, $\theta$ is a finite supremum of seminorms, hence a
seminorm, convex and continuous.  For unitaries $U,V$ the bijection $(X,Y)\mapsto(V^*XU^*,V^*YU^*)$ of $M_n(\C)^2$
preserves the Frobenius norms, and
\[
T_{UZV}(V^*XU^*,V^*YU^*)=V^*T_Z(X,Y)U^*,\qquad \mathsf W_{UZV}(V^*XU^*,V^*YU^*)=V^*\mathsf W_Z(X,Y)U^*
\]
(for the second, $(UZV)^*=V^*Z^*U^*$).  The map $C\mapsto V^*CU^*$ is unitary on $M_n(\C)$ and sends $Y$ to
$V^*YU^*$, so it carries $P^\perp_Y\mathsf W_Z(X,Y)$ to $P^\perp_{V^*YU^*}\mathsf W_{UZV}(V^*XU^*,V^*YU^*)$.  Hence
$\mathcal L_{UZV}(V^*XU^*,V^*YU^*)=\mathcal L_Z(X,Y)$, and $\theta(UZV)=\theta(Z)$.  By the singular value decomposition,
$\theta(Z)=\theta(\diag(s_1(Z),\dots,s_n(Z)))$, and $\theta(\diag(t))$, $t\in\R^n$, is unchanged under permutations of the
entries of $t$ ($U$ a permutation matrix, $V=U^*$) and changes of their signs ($U=I$, $V$ diagonal with entries
$\pm1$).  For $Z=\diag(s_1(Z),\dots,s_n(Z))$, $X=E_{12}$ and $Y=E_{21}$ give $T_Z(X,Y)=s_2(Z)E_{11}-s_1(Z)E_{22}$, so
$\mathfrak n(Z)\ge\norm Z_{(2),2}$, for every $Z$ by invariance; and $\phi\ge\mathfrak n$ since
$\mathcal L_Z(X,Y)\ge\norm{T_Z(X,Y)}^2$.
\item \emph{Averaging over signs.}  Let $Z\ne0$ and $s_k=s_k(Z)$, so that $\theta(Z)=\theta(\diag(s))$ by step~1.  Let
$0<\beta\le s_1$ with $\beta\ge s_2$ (for instance $\beta=s_2$ if $s_2>0$).  Lemma~\ref{lem:red-sign} writes
$\diag(s)=\sum_\epsilon w_\epsilon D_\epsilon$ with $D_\epsilon=\diag(s_1,\beta,\dots,\beta)E_\epsilon$, $E_\epsilon$ a diagonal
sign matrix, so $\theta(D_\epsilon)=\theta(\diag(s_1,\beta,\dots,\beta))$ by step~1, and this is at most
$(s_1^2+\beta^2)^{1/2}$ by the hypothesis.  Since $\theta$ is a seminorm,
\[
\theta(Z)=\theta\Bigl(\sum_\epsilon w_\epsilon D_\epsilon\Bigr)\le\sum_\epsilon w_\epsilon\,\theta(D_\epsilon)\le(s_1^2+\beta^2)^{1/2}.
\]
If $s_2>0$, take $\beta=s_2$: then $\theta(Z)\le\norm Z_{(2),2}$.  If $s_2=0$, let $\beta\downarrow0$: again
$\theta(Z)\le s_1=\norm Z_{(2),2}$.  With step~1, and trivially for $Z=0$, $\theta(Z)=\norm Z_{(2),2}$ for every $Z$.\qedhere
\end{enumerate}
\end{proof}
Geometrically, step~2 is the extreme point principle for the convex function $\theta$: up to order and signs,
the vectors $(a,b,\dots,b)$ with $a\ge b\ge0$ and $a^2+b^2=1$ are the extreme points of the unit ball
$\{t\in\R^n:\max_{i<j}(t_i^2+t_j^2)\le1\}$ of the Ky Fan $(2),2$ norm on diagonal matrices, and
Lemma~\ref{lem:red-sign} writes a diagonal matrix explicitly as a convex combination of matrices with
singular values $(s_1,\beta,\dots,\beta)$.  The proof uses only this explicit decomposition, which also
transfers the cases of equality (Lemma~\ref{lem:eq-transfer}).

\begin{remark}[the adjoint in $\mathsf W_Z$]\label{rem:red-two-level-2V}
For Hermitian $Z$ the matrix $\mathsf W_Z$ is the matrix $\mathsf W$ of Theorem~\ref{thm:intro-2V} with $Z$ in place
of $D$; for other $Z$ the adjoint matters.  For $n=2$, $Z=\diag(1,i)$, $X=E_{21}$ and $Y=Z$ one has
$T_Z(X,Y)=2E_{21}$ and $ZX^*Y+YX^*Z=2i\,E_{12}\perp Y$, so the left side of \eqref{eq:intro-2V} written literally
with $Z$ in place of $D$ equals $8=4\norm X^2\norm Y^2$, whereas $\norm Z_{(2),2}^2\norm X^2\norm Y^2=4$; here
$\mathsf W_Z(X,Y)=Z^*E_{12}Z+ZE_{12}Z^*=0$.  Thus \eqref{eq:intro-2V} with $Z$ in place of $D$ fails for this normal
$Z$, while with $\mathsf W_Z$ it holds for every $Z$ (Corollary~\ref{cor:m2-2V-all}).
\end{remark}

For weights with $\lambda_2=\dots=\lambda_n$, (2V) localizes to a subspace of dimension at most three.  For a
subspace $V$ let $\Pi_V$ be the orthogonal projection onto $V$, and for $X\in M_n(\C)$ write $\Upsilon_X$ for the
functional $\Upsilon$ of Section~\ref{subsec:red-2V} formed with $X$.

\begin{lemma}[localization]\label{lem:red-local}
Let $D=bI+(a-b)ee^*$ with $a\ge b>0$ and a unit vector $e\in\C^n$, and $c=a^2+b^2$.  Let $\psi_1,\psi_2$ be
unit vectors, let $V$ be a subspace containing $e,\psi_1,\psi_2$, $\Pi=\Pi_V$ and $\Pi^\perp=I-\Pi$, and for
$X\in M_n(\C)$ put $X_{VV}=\Pi X\Pi$, $X_{VW}=\Pi X\Pi^\perp$, $X_{WV}=\Pi^\perp X\Pi$ and
$X_{WW}=\Pi^\perp X\Pi^\perp$.
\begin{enumerate}
\item[(i)] $\Upsilon_X(\psi_1,\psi_2)=\Upsilon_{X_{VV}}(\psi_1,\psi_2)+\Xi_{VW}+\Xi_{WV}$, where
\[
\Xi_{VW}=\norm{X_{VW}^*D\psi_2}^2+b^2\norm{X_{VW}^*\psi_1}^2,\qquad
\Xi_{WV}=\norm{X_{WV}D\psi_1}^2+b^2\norm{X_{WV}\psi_2}^2 ;
\]
in particular $X_{WW}$ does not occur.
\item[(ii)] $\Xi_{VW}\le c\norm{X_{VW}}^2$ and $\Xi_{WV}\le c\norm{X_{WV}}^2$.  If $a>b$ and $\psi_1,\psi_2$ are not
both in $\C e$, then $\Xi_{VW}=c\norm{X_{VW}}^2$ only if $X_{VW}=0$, and $\Xi_{WV}=c\norm{X_{WV}}^2$ only if
$X_{WV}=0$.
\end{enumerate}
\end{lemma}

\begin{proof}
Since $e\in V$, $Dv=bv+(a-b)\ip ev e\in V$ for $v\in V$, and $V^\perp\subseteq e^\perp$; hence $D\Pi=\Pi D$ and
$D\Pi^\perp=b\Pi^\perp$.  Also $\Pi\psi_k=\psi_k$.

(i) We use Lemma~\ref{lem:red-decomp} for $X$ and for $X_{VV}$.  Put $Y=\psi_1\psi_2^*$.  Since $D\psi_1,D\psi_2\in V$,
we have $XD\psi_1=X_{VV}D\psi_1+X_{WV}D\psi_1$ and $X^*D\psi_2=X_{VV}^*D\psi_2+X_{VW}^*D\psi_2$, so
\begin{align*}
T_D(X,Y)&=XD\psi_1\psi_2^*-\psi_1(X^*D\psi_2)^*\\
&=T_D(X_{VV},Y)+X_{WV}D\psi_1\psi_2^*-\psi_1(X_{VW}^*D\psi_2)^* .
\end{align*}
The three terms lie in $\Pi M_n\Pi$, $\Pi^\perp M_n\Pi$ and $\Pi M_n\Pi^\perp$, which are mutually orthogonal, so
\[
\norm{T_D(X,Y)}^2=\norm{T_D(X_{VV},Y)}^2+\norm{X_{WV}D\psi_1}^2+\norm{X_{VW}^*D\psi_2}^2 .
\]
Next,
$X^*\psi_1=X_{VV}^*\psi_1+X_{VW}^*\psi_1$ with $X_{VW}^*\psi_1\in V^\perp$, so
$DX^*\psi_1=DX_{VV}^*\psi_1+bX_{VW}^*\psi_1$.  The first term lies in $V$, and so does its image under
$P^\perp_{\psi_1}$ (as $\psi_1\in V$); the second is orthogonal to $V\ni\psi_1$ and is fixed by $P^\perp_{\psi_1}$.  Hence
the first of the equations
\begin{align*}
\norm{P^\perp_{\psi_1}DX^*\psi_1}^2&=\norm{P^\perp_{\psi_1}DX_{VV}^*\psi_1}^2+b^2\norm{X_{VW}^*\psi_1}^2,\\
\norm{P^\perp_{\psi_2}DX\psi_2}^2&=\norm{P^\perp_{\psi_2}DX_{VV}\psi_2}^2+b^2\norm{X_{WV}\psi_2}^2
\end{align*}
holds, and the second follows in the same way from $X\psi_2=X_{VV}\psi_2+X_{WV}\psi_2$.  Adding the three
equations and using Lemma~\ref{lem:red-decomp} for $X$ and for $X_{VV}$ gives (i).

(ii) Since $\norm D_{\mathrm{op}}=a$ and the operator norm is at most the Frobenius norm,
$\norm{X_{VW}^*D\psi_2}\le\norm{X_{VW}}\norm{D\psi_2}\le a\norm{X_{VW}}$ and $\norm{X_{VW}^*\psi_1}\le\norm{X_{VW}}$; so
$\Xi_{VW}\le(a^2+b^2)\norm{X_{VW}}^2$.  Suppose that equality holds and $X_{VW}\ne0$.  Since $b>0$, both estimates are
equalities.  From $\norm{X_{VW}}^2-\norm{X_{VW}^*\psi_1}^2=\norm{P^\perp_{\psi_1}X_{VW}}^2=0$ we get $X_{VW}=\psi_1u^*$ with
$u=X_{VW}^*\psi_1\ne0$; then $X_{VW}^*D\psi_2=\ip{\psi_1}{D\psi_2}u$, and the first equality gives
$\abs{\ip{\psi_1}{D\psi_2}}=a$.  As $\abs{\ip{\psi_1}{D\psi_2}}\le\norm{D\psi_2}\le a$ and
$\norm{D\psi_2}^2=b^2+(a^2-b^2)\abs{\ip e{\psi_2}}^2$, this forces $\psi_2\in\C e$ (because $a>b$) and then
$\psi_1\in\C D\psi_2=\C e$, which is excluded.  The claims for $\Xi_{WV}$ are the same claims for $X_{WV}^*$ in place of
$X_{VW}$, with $\psi_1$ and $\psi_2$ exchanged.
\end{proof}

\begin{corollary}[reduction to $n=3$]\label{cor:red-local}
Let $a\ge b>0$.  If (2V) holds for $\diag(a,b,b)\in M_3(\C)$, then it holds for $\diag(a,b,\dots,b)\in M_n(\C)$
for every $n\ge2$.
\end{corollary}

\begin{proof}
Two facts follow directly from the definition of $\Upsilon$ (through $\mathsf M_1$, $\mathsf M_2$ and $\chi$).  If $U$ is
unitary and $U^*DU=D'$, then $\Upsilon$ formed with $D',U^*XU$ at $(U^*\psi_1,U^*\psi_2)$ equals $\Upsilon$ formed with
$D,X$ at $(\psi_1,\psi_2)$.  If $D=D_1\oplus D_2$, $X=X_1\oplus0$ and $\psi_k=\psi_k'\oplus0$, then $\Upsilon$ equals the
functional formed with $D_1,X_1$ at $(\psi_1',\psi_2')$.  In both cases $c$ and $\norm X$ do not change.

Let $D=\diag(a,b,\dots,b)=bI+(a-b)ee^*$ with $e=\mathbf e_1$, the first standard basis vector of $\C^n$, and
$c=a^2+b^2$.  If $n=2$, the second fact with
$\diag(a,b,b)=D\oplus b$, $X\oplus0$ and $\psi_k\oplus0$ gives (2V) for $D$.  Let $n\ge3$, and fix $X$ and unit vectors
$\psi_1,\psi_2$.  Choose a three-dimensional subspace $V$ containing $e,\psi_1,\psi_2$, an orthonormal basis
$(e,v_2,v_3)$ of $V$ and an orthonormal basis of $V^\perp$, and let $U$ be the unitary matrix with these $n$ vectors as
columns.  Since $Dv_k=bv_k$ and $D=b$ on $V^\perp$, $U^*DU=D=\diag(a,b,b)\oplus bI_{n-3}$, and $U^*X_{VV}U=X'\oplus0$,
$U^*\psi_k=\psi_k'\oplus0$ with $X'\in M_3(\C)$, $\psi_k'\in\C^3$.  By the two facts and (2V) for $\diag(a,b,b)$,
$\Upsilon_{X_{VV}}(\psi_1,\psi_2)\le c\norm{X'}^2=c\norm{X_{VV}}^2$.  Lemma~\ref{lem:red-local} gives
\[
\Upsilon_X(\psi_1,\psi_2)\le c\bigl(\norm{X_{VV}}^2+\norm{X_{VW}}^2+\norm{X_{WV}}^2\bigr)\le c\norm X^2 .\qedhere
\]
\end{proof}

The localization uses that every subspace containing $e$ is invariant under $D=bI+(a-b)ee^*$ and that $D$ is
scalar on its orthogonal complement.  For weights with more than two distinct eigenvalues no subspace of fixed
dimension containing $\psi_1,\psi_2$ need have these properties; this is why the reduction to
$\lambda_2=\dots=\lambda_n$ (Proposition~\ref{prop:red-two-level}) comes first.  The plane of the two least weighted
levels, which carries the extremal pairs (Theorem~\ref{thm:intro-eq}), enters the proof only through $c=d_1^2+d_2^2$
and the Ky Fan maximum principle (Lemma~\ref{lem:red-KyFan} below), which produces the constraint set $K^\C_2$ of
Definition~\ref{def:F}.

\subsection{Adapted bases}\label{subsec:red-adapted}

\emph{From now on until the end of Section~\ref{subsec:red-assembly}, $n=3$}, so that
$D=\diag(d_1\ge d_2\ge d_3>0)$; Section~\ref{subsec:red-two-level} has reduced the weights we need to this case.

\emph{Normalisation.}  $\mathsf M_1$ and $\mathsf M_2$ are quadratic and $\chi$ is linear in $D$, so, writing the weight
as a subscript (from here on a subscript on $\Upsilon$ always denotes the weight; the notation $\Upsilon_X$ of
Section~\ref{subsec:red-two-level} is not used again), for $\Lambda>0$
\begin{equation}\label{eq:red-normal}
\Upsilon_{D/\sqrt\Lambda}(\psi_1,\psi_2)=\Lambda^{-1}\,\Upsilon_D(\psi_1,\psi_2) .
\end{equation}
Thus (2V)$[\Lambda]$ for $D$ is (2V)$[1]$ for $D/\sqrt\Lambda$, whose constant $c$ is divided by $\Lambda$.  Until the end of
the Ky Fan lemma (Lemma~\ref{lem:red-KyFan}) we study (2V)$[1]$, that is, $\Upsilon(\psi_1,\psi_2)\le\norm X^2$; the
number $\Lambda$ returns only in Section~\ref{subsec:red-assembly}.  Replacing $D$ by $D/\sqrt\Lambda$ does not change the
adapted bases below and divides $h$, $r$, $s$ and $B$ by $\sqrt\Lambda$.

\begin{definition}[adapted bases \cite{Nicolas2026}]\label{def:red-adapted}
For unit vectors $\psi_1,\psi_2\in\C^3$ put
\[
h=\langle\psi_1,D\psi_2\rangle\in\C,\qquad r=\norm{P^\perp_{\psi_1}D\psi_2}\ge0,\qquad
s=\norm{P^\perp_{\psi_2}D\psi_1}\ge0 .
\]
Let $e_1=\psi_1$; if $r>0$ let $e_2=P^\perp_{\psi_1}D\psi_2/r$, and if $r=0$ let $e_2$ be any unit
vector orthogonal to $\psi_1$; let $e_3$ be a unit vector orthogonal to $e_1,e_2$.
Likewise let $f_1=\psi_2$, $f_2=P^\perp_{\psi_2}D\psi_1/s$ if $s>0$ (any unit vector orthogonal to
$\psi_2$ if $s=0$), and $f_3\perp f_1,f_2$ a unit vector.  Let
$\mathsf E=(e_1,e_2,e_3)$ and $\mathsf F=(f_1,f_2,f_3)$ (unitary), $\mathbf G=\mathsf E^*D\mathsf F$
and $\mathsf Z=\mathsf E^*X\mathsf F$, so that $\mathbf G_{ij}=\langle e_i,Df_j\rangle$,
$z_{ij}=\langle e_i,Xf_j\rangle$ and $\norm{\mathsf Z}=\norm X$.  Write
$p=(z_{21},z_{31})^\T$ and $q=(z_{12},z_{13})^\T$ in $\C^2$.  Outside such index formulas we write $\psi_1,\psi_2$
for $e_1,f_1$; in this section $e_1$ denotes the first standard basis vector of $\C^2$, and $\mathbf e_1$ that
of $\C^3$.
\end{definition}

\begin{lemma}[zero pattern]\label{lem:red-zero}
In an adapted basis
\[
\mathbf G=\begin{pmatrix}h&s\,e_1^\T\\ r\,e_1&B\end{pmatrix}\qquad(B\in M_2(\C)),
\]
and $\mathbf G\mathbf G^*=\mathsf E^*D^2\mathsf E$, $\mathbf G^*\mathbf G=\mathsf F^*D^2\mathsf F$ have
eigenvalues $d_1^2,d_2^2,d_3^2$.  Here $h$ and $B$ are complex and $r,s$ are nonnegative reals.
\end{lemma}

\begin{proof}
$D\psi_2=h\psi_1+re_2$ (also when $r=0$) gives the first column.  Since $D$ is Hermitian,
$\langle\psi_2,D\psi_1\rangle=\bar h$ and $D\psi_1=\bar h\psi_2+sf_2$, so
$\mathbf G_{1j}=\overline{\langle f_j,D\psi_1\rangle}$ equals $h,s,0$.  Finally
$\mathbf G\mathbf G^*=\mathsf E^*D^2\mathsf E$ and $\mathbf G^*\mathbf G=\mathsf F^*D^2\mathsf F$.
\end{proof}
The numbers $h,r,s$, the adapted bases and the first row and column of $\mathbf G$ are those of
\cite[Setup~2.1]{Nicolas2026} (with $u=\psi_1$, $v=\psi_2$), where they are used for a rank-one factor.

The phases of $e_2,f_2$ were used to make $r,s$ real; the phase of $h$ cannot be chosen, but $h$
enters below only through $\abs h^2$ and inside completed squares.

\begin{lemma}[$\Upsilon$ in an adapted basis]\label{lem:red-coord}
\begin{multline*}
\Upsilon(\psi_1,\psi_2)=\abs{sq_1-rp_1}^2+\sum_{i=2,3}\abs{\bar hz_{i1}+sz_{i2}}^2
+\sum_{j=2,3}\abs{\bar hz_{1j}+rz_{2j}}^2\\
+\norm{r\bar z_{11}e_1+B\bar q}^2+\norm{sz_{11}e_1+B^*p}^2 .
\end{multline*}
\end{lemma}

\begin{proof}
We write the three terms of Lemma~\ref{lem:red-decomp} in coordinates.

\emph{The commutator.}
$T_D(X,\psi_1\psi_2^*)=\mathbf c_1\psi_2^*-\psi_1\mathbf c_2^*$ with $\mathbf c_1=XD\psi_1$, $\mathbf c_2=X^*D\psi_2$.
From $D\psi_1=\bar h\psi_2+sf_2$ and $D\psi_2=h\psi_1+re_2$ we get $\langle e_i,\mathbf c_1\rangle=\bar hz_{i1}+sz_{i2}$
and $\langle\mathbf c_2,f_j\rangle=\bar hz_{1j}+rz_{2j}$.  So the nonzero entries of
$\mathsf E^*T_D(X,\psi_1\psi_2^*)\mathsf F$ lie in the first column ($i\ge2$), the first row ($j\ge2$) and at
$(1,1)$, where the entry is $(\bar hz_{11}+sz_{12})-(\bar hz_{11}+rz_{21})=sq_1-rp_1$.

\emph{The second remainder.}  $X^*\psi_1=\sum_j\overline{z_{1j}}f_j$ and $Df_j=\sum_i\mathbf G_{ij}e_i$,
so $DX^*\psi_1=\sum_i(\mathbf G\bar{\mathbf a})_ie_i$ with $\mathbf a=(z_{11};q)$.  The projection
$P^\perp_{\psi_1}$ removes the $\psi_1$ component, and the rows of $\mathbf G$ with index $\ge2$ are
$(r\,e_1\ \ B)$, so what remains is $r\bar z_{11}e_1+B\bar q$.

\emph{The third remainder.}  $DX\psi_2=\sum_iz_{i1}De_i$ and
$De_i=\sum_j\overline{\mathbf G_{ij}}f_j$, so the coordinates of $DX\psi_2$ in the basis $(f_j)$ are
$\mathbf G^*\boldsymbol\xi$ with $\boldsymbol\xi=(z_{11};p)$.  The columns of $\mathbf G$ with index
$\ge2$ are $(s\,e_1^\T;B)$, so after removing the $\psi_2$ component what remains is
$sz_{11}e_1+B^*p$.
\end{proof}

Only the fourth term contains $\bar z_{11}$ and $\bar q$: the second remainder comes from $DX^*\psi_1$
and is conjugate linear in the entries of $\mathsf Z$, the third comes from $DX\psi_2$ and is linear.  The
right-hand side involves only the first two rows and the first two columns of $\mathsf Z$.

\subsection{Completing squares}

\begin{lemma}[elimination of one variable]\label{lem:red-elim}
Let $\mathbf A\in\C^k$, $\ell\in\R^k$, $\Lambda>\norm\ell^2$ and $t\in\C$.  Then
\[
\norm{\mathbf A+t\ell}^2-\Lambda\abs t^2=\norm{\mathbf A}^2
+\frac{\abs{\ell^\T\mathbf A}^2}{\Lambda-\norm\ell^2}-(\Lambda-\norm\ell^2)\abs{t-t^\star}^2,
\qquad t^\star=\frac{\ell^\T\mathbf A}{\Lambda-\norm\ell^2}.
\]
\end{lemma}

\begin{proof}
Since $\ell$ is real, $\langle\mathbf A,t\ell\rangle=t\,\overline{\ell^\T\mathbf A}$, so
$\norm{\mathbf A+t\ell}^2=\norm{\mathbf A}^2+2\operatorname{Re}(\bar t\,\ell^\T\mathbf A)+\abs t^2\norm\ell^2$.
Expanding the last term on the right gives
$-(\Lambda-\norm\ell^2)\abs t^2+2\operatorname{Re}(\bar t\,\ell^\T\mathbf A)
-\frac{\abs{\ell^\T\mathbf A}^2}{\Lambda-\norm\ell^2}$, and the two sides agree.
\end{proof}

We first define, for every $m\ge1$, the parameters and matrices that occur.  This paper uses them only for $m=2$;
the case $m\ge3$ is treated in the separate paper (in preparation) mentioned in Section~\ref{subsec:intro-second}.

\begin{definition}[parameters]\label{def:params}
Let $R,S\ge0$, $H\ge0$ with $G:=1-R-S-H>0$, and put $r=\sqrt R$, $s=\sqrt S$,
\[
\nu=G+H,\quad\sigma=1-S,\quad\varrho=1-R,\quad a_1=G+S,\quad a_2=G+R,\quad
\gamma_p=\frac{a_2}\sigma,\quad\gamma_q=\frac{a_1}\varrho .
\]
All of these numbers except $r,s$ are positive ($r=0$ or $s=0$ is allowed).  We also put
\[
a_{12}=a_1+a_2=1+G-H>0,\qquad c_D=2\nu G+\nu R+SG>0 .
\]
\end{definition}

\begin{definition}[the matrices $F$ and $F_\C$]\label{def:F}\label{def:FC}
Let $m\ge1$, let $e_1$ be the first standard basis vector of $\C^m$, and put $b=(se_1;\,re_1)\in\R^{2m}$,
$D_p=\diag(G,a_2,\dots,a_2)$ and $D_q=\diag(G,a_1,\dots,a_1)$ (all $m\times m$).  For $B\in M_m(\C)$ let
$\mathbf m(B)=(sBe_1;\,rB^\T e_1)\in\C^{2m}$ (a transpose, not an adjoint),
\begin{gather*}
\widehat S_p=\diag(G,\gamma_p,\dots,\gamma_p)-BB^*,\qquad \widehat S_q=\diag(G,\gamma_q,\dots,\gamma_q)-B^\T\bar B,\\
F_\C(B)=\begin{pmatrix}\widehat S_p&0\\0&\widehat S_q\end{pmatrix}
+\frac1\nu\bigl(G\,bb^\T-\mathbf m(B)\mathbf m(B)^*\bigr),
\end{gather*}
$S_p=D_p-BB^*$, $S_q=D_q-B^*B$ and $K_m^\C=\{B\in M_m(\C):S_p\succeq0,\ S_q\succeq0\}$.  For real $B$
we write $F(B)=F_\C(B)$, which is then real symmetric ($B^\T\bar B=B^\T B$), and $K_m=K_m^\C\cap M_m(\R)$.
Equivalently, with $\Lambda_p=\diag(\frac G\nu,\gamma_p,\dots,\gamma_p)$, $\Lambda_q=\diag(\frac G\nu,\gamma_q,\dots,\gamma_q)$
and $\hat w=(re_1;\,-se_1)\in\R^{2m}$,
\[
F_\C(B)=\begin{pmatrix}\Lambda_p-BB^*&0\\0&\Lambda_q-B^\T\bar B\end{pmatrix}
-\frac1\nu\bigl(\mathbf m(B)\mathbf m(B)^*+G\,\hat w\hat w^\T\bigr),
\]
because $\frac G\nu-G=\frac G\nu(R+S)$ and, on the two coordinates $(e_1;0)$, $(0;e_1)$, the Lagrange identity
$(R+S)I_2-\hat w_E\hat w_E^\T=b_Eb_E^\T$ holds with $\hat w_E=(r,-s)^\T$, $b_E=(s,r)^\T$.
\end{definition}

\begin{theorem}[completing squares]\label{thm:red-square}
Let $n=3$, let $\psi_1,\psi_2$ be unit vectors, and let $h,r,s,B$ and $\mathsf Z$ be as in
Definition~\ref{def:red-adapted} and Lemma~\ref{lem:red-zero}.  Assume $r^2+s^2<1$, put $R=r^2$, $S=s^2$, $H=\abs h^2$, and use the notation of Definitions~\ref{def:params}
and~\ref{def:F} with $m=2$; here $G=1-R-S-H$ may be $\le0$, and only $\nu,\sigma,\varrho>0$ are needed.  With
$\mathbf x=(p;q)$, for every $\mathsf Z$
\begin{multline}\label{eq:red-square}
\norm{\mathsf Z}^2-\Upsilon(\psi_1,\psi_2)=\abs{z_{33}}^2+\mathbf x^*F_\C(B)\,\mathbf x
+\sigma\Bigl\lvert z_{32}-\frac{s\bar hz_{31}}\sigma\Bigr\rvert^2\\
+\varrho\Bigl\lvert z_{23}-\frac{r\bar hz_{13}}\varrho\Bigr\rvert^2
+\nu\Bigl\lvert z_{22}-\frac{\bar h(sp_1+rq_1)}\nu\Bigr\rvert^2
+\nu\abs{z_{11}-t^\star_{11}}^2,
\end{multline}
where $t^\star_{11}=\bigl(r(\bar Bq)_1+s(B^*p)_1\bigr)/\nu=\mathbf m(B)^*\mathbf x/\nu$.  Consequently, if
$F_\C(B)\succeq0$, then \textup{(2V)}$[1]$ holds for the pair $(\psi_1,\psi_2)$.
\end{theorem}

\begin{proof}
Split $\norm{\mathsf Z}^2=\abs{z_{11}}^2+\norm p^2+\norm q^2+\abs{z_{22}}^2+\abs{z_{32}}^2
+\abs{z_{23}}^2+\abs{z_{33}}^2$, and sort the terms of
Lemma~\ref{lem:red-coord} by the variables they contain.
\begin{enumerate}
\item The arm entry $z_{32}$ occurs only in $\abs{\bar hp_2+sz_{32}}^2$ ($z_{31}=p_2$).
Lemma~\ref{lem:red-elim} with $k=1$, $\Lambda=1$, $\mathbf A=\bar hp_2$, $\ell=s$ gives
$\abs{z_{32}}^2-\abs{\bar hp_2+sz_{32}}^2
=-\frac H\sigma\abs{p_2}^2+\sigma\abs{z_{32}-t^\star}^2$.  The arm entry $z_{23}$ gives
$-\frac H\varrho\abs{q_2}^2$ plus a square in the same way.
\item $z_{22}$ occurs only in $\abs{\bar hp_1+sz_{22}}^2$ (the term $i=2$) and
$\abs{\bar hq_1+rz_{22}}^2$ (the term $j=2$).  Lemma~\ref{lem:red-elim} with $k=2$, $\Lambda=1$,
$\mathbf A=(\bar hp_1;\bar hq_1)$, $\ell=(s;r)$, $\ell^\T\mathbf A=\bar h(sp_1+rq_1)=\bar h\,b^\T\mathbf x$ (condition
$R+S<1$) gives the value
$-H\bigl(\abs{p_1}^2+\abs{q_1}^2\bigr)-\frac H\nu\abs{b^\T\mathbf x}^2$ plus a
square.
\item $z_{11}$ occurs only in the last two norms.  We conjugate the first one:
$\norm{r\bar z_{11}e_1+B\bar q}=\norm{rz_{11}e_1+\bar Bq}$ ($r$ is real).  Now both norms are
complex linear in $z_{11}$, and Lemma~\ref{lem:red-elim} with $k=4$, $\Lambda=1$,
$\mathbf A=(\bar Bq;B^*p)$, $\ell=(r\,e_1;s\,e_1)$, $\norm\ell^2=R+S$,
$\ell^\T\mathbf A=r(\bar Bq)_1+s(B^*p)_1=\mathbf m(B)^*\mathbf x$ gives the value
$-\norm{\bar Bq}^2-\norm{B^*p}^2-\frac1\nu\abs{\mathbf m(B)^*\mathbf x}^2$ plus a square.  Since
$(\bar B)^*=B^\T$, this value is
$-p^*BB^*p-q^*B^\T\bar Bq-\frac1\nu\abs{\mathbf m(B)^*\mathbf x}^2$.
\item $\abs{sq_1-rp_1}^2$ involves only $p,q$, and $z_{33}$ does not occur in $\Upsilon$.
\end{enumerate}
The groups of variables in 1--4 occur in disjoint sets of terms, so $\norm{\mathsf Z}^2-\Upsilon$ is
the sum of $\abs{z_{33}}^2$, the four squares, and the Hermitian form $Q(\mathbf x)$ formed by
$\norm{\mathbf x}^2=\norm p^2+\norm q^2$, the values in 1--3 and $-\abs{sq_1-rp_1}^2$.

It remains to show $Q(\mathbf x)=\mathbf x^*F_\C(B)\mathbf x$.  By the Lagrange identity
$\abs{sq_1-rp_1}^2=(R+S)\bigl(\abs{p_1}^2+\abs{q_1}^2\bigr)-\abs{sp_1+rq_1}^2$ (the cross terms
$\mp2rs\operatorname{Re}(q_1\bar p_1)$ cancel), the terms of $Q$ in $p_1,q_1$ alone are
\[
(1-H-R-S)\bigl(\abs{p_1}^2+\abs{q_1}^2\bigr)+\Bigl(1-\frac H\nu\Bigr)\abs{b^\T\mathbf x}^2
=G\bigl(\abs{p_1}^2+\abs{q_1}^2\bigr)+\frac G\nu\,\mathbf x^*bb^\T\mathbf x ,
\]
since $\nu=G+H$ and $b$ is real.  The coefficient of $\abs{p_2}^2$ is $1-\frac H\sigma=\frac{a_2}\sigma=\gamma_p$,
and that of $\abs{q_2}^2$ is $\gamma_q$.  With step~3,
\begin{multline*}
Q(\mathbf x)=\mathbf x^*\Bigl[\bigl(\diag(G,\gamma_p)-BB^*\bigr)\oplus\bigl(\diag(G,\gamma_q)-B^\T\bar B\bigr)\\
+\frac1\nu\bigl(Gbb^\T-\mathbf m(B)\mathbf m(B)^*\bigr)\Bigr]\mathbf x=\mathbf x^*F_\C(B)\mathbf x ,
\end{multline*}
which proves \eqref{eq:red-square}.  If $F_\C(B)\succeq0$, every term on the right of
\eqref{eq:red-square} is nonnegative ($\sigma,\varrho,\nu>0$), so $\Upsilon(\psi_1,\psi_2)\le\norm{\mathsf Z}^2=\norm X^2$.
\end{proof}

Conversely, taking $z_{33}=0$ and the eliminated variables at their optimal values shows that (2V)$[1]$ for the
pair $(\psi_1,\psi_2)$ implies $F_\C(B)\succeq0$, since $p$ and $q$ are arbitrary; only the implication of
Theorem~\ref{thm:red-square} is used below.

\begin{remark}[orientation and the transpose]\label{rem:red-orient}
The transposes in $\mathbf m(B)=(sBe_1;rB^\T e_1)$ and in $\widehat S_q$ come from the
conjugation of one remainder in step~3 of the proof of Theorem~\ref{thm:red-square}, and they cannot be replaced
by adjoints: for $m=1$, $R=S=\frac14$, $H=0$ and $B=(i/\sqrt2)\in K^\C_1$ one has $F_\C(B)=0$, while the form with
$B^*B$ in place of $B^\T\bar B$ and border $(sBe_1;rB^*e_1)$ is
$\left(\begin{smallmatrix}0&1/2\\1/2&0\end{smallmatrix}\right)$, with eigenvalue $-\frac12$.  The set $K^\C_m$ does
not depend on the choice, since $B^\T\bar B=\overline{B^*B}$.
\end{remark}

\subsection{The Ky Fan maximum principle and the constraint set}

Put $\tau_\psi=\abs h^2+r^2+s^2$.

\begin{lemma}[Ky Fan]\label{lem:red-KyFan}
In an adapted basis:
\begin{enumerate}
\item[(a)] $BB^*+r^2E_{11}\preceq(c-\abs h^2-s^2)I$ and $B^*B+s^2E_{11}\preceq(c-\abs h^2-r^2)I$;
\item[(b)] $\tau_\psi\le c$;
\item[(c)] if $c<1$, then, with $R=r^2$, $S=s^2$, $H=\abs h^2$ in Definition~\ref{def:F},
$S_p=D_p-BB^*\succeq(1-c)I\succ0$ and $S_q=D_q-B^*B\succeq(1-c)I\succ0$;
in particular $B\in K^\C_2$, and moreover $G=1-\tau_\psi\ge1-c>0$,
$\nu=G+H>0$ and $r^2+s^2<1$.
\end{enumerate}
\end{lemma}

\begin{proof}
The Ky Fan maximum principle \cite[Problem~I.6.15, Exercise~II.1.13]{BhatiaMA} says: if $N$ is a Hermitian
matrix, then for orthonormal $x,y$, $x^*Nx+y^*Ny=\Tr(PN)$ is at most the sum of
the two largest eigenvalues of $N$, where $P$ is the projection onto $\operatorname{span}(x,y)$.  For
$N=\mathbf G\mathbf G^*$ this sum is $d_1^2+d_2^2=c$ (Lemma~\ref{lem:red-zero}), and
$x^*\mathbf G\mathbf G^*x=\norm{\mathbf G^*x}^2$.

(a) Take $x=\mathbf e_1$ and $y=(0;v)$ with $v\in\C^2$ a unit vector.  Since
$\mathbf G^*=\begin{pmatrix}\bar h&r\,e_1^\T\\ s\,e_1&B^*\end{pmatrix}$, we get
$\norm{\mathbf G^*\mathbf e_1}^2=\abs h^2+s^2$ and $\mathbf G^*(0;v)=(rv_1;B^*v)$, so
$\norm{\mathbf G^*(0;v)}^2=v^*(r^2E_{11}+BB^*)v$.  Ky Fan gives
$v^*(BB^*+r^2E_{11})v\le c-\abs h^2-s^2$ for every unit vector $v$, which is the first operator
inequality.  For the second use $N=\mathbf G^*\mathbf G$ (same eigenvalues):
$\norm{\mathbf G\mathbf e_1}^2=\abs h^2+r^2$, $\mathbf G(0;v)=(sv_1;Bv)$ and
$\norm{\mathbf G(0;v)}^2=v^*(s^2E_{11}+B^*B)v$.

(b) Put $v=e_1$ in the first inequality of (a):
$\norm{B^*e_1}^2+r^2\le c-\abs h^2-s^2$, so $\tau_\psi\le c-\norm{B^*e_1}^2\le c$.  (The bound
$\tau_\psi\le c$ is \cite[Theorem~3.1]{Nicolas2026}; there $\tau_\psi$ appears as the squared norm of the
twisted commutator map with the rank-one factor $\psi_1\psi_2^*$ \cite[Theorem~2.2]{Nicolas2026}, and it is
bounded by the same Ky Fan argument.)

(c) The first inequality of (a) reads $BB^*+RE_{11}\preceq(c-H-S)I$.  Since
$D_p=(1-S-H)I-RE_{11}$ (Definition~\ref{def:F}), $S_p=D_p-BB^*\succeq(1-c)I$.  The
same argument with the second inequality and $D_q=(1-R-H)I-SE_{11}$ gives $S_q$.  The rest follows
from (b): $G=1-\tau_\psi$ and $1>c\ge\tau_\psi\ge r^2+s^2$.
\end{proof}

\subsection{Assembly and the limit \texorpdfstring{$\Lambda\downarrow c$}{Lambda to c}}\label{subsec:red-assembly}

\begin{theorem}[reduction]\label{thm:red-chain}
Let $n=3$, and suppose that $F_\C(B)\succeq0$ for all parameters as in Definition~\ref{def:params}
\textup{(}in particular $G>0$\textup{)} and all $B\in K^\C_2$.  Then (2V) holds for every
$D=\diag(d_1\ge d_2\ge d_3>0)$.
\end{theorem}
The hypothesis is Theorem~\ref{thm:m2-D}~(iii), proved in Section~\ref{sec:m2}, whose proof does not use
Theorem~\ref{thm:red-chain}.

\begin{proof}
Fix $X$ and unit vectors $\psi_1,\psi_2$, and take adapted bases (Definition~\ref{def:red-adapted}).
\begin{enumerate}
\item Let $\Lambda>c$ be arbitrary.  The weight $D'=D/\sqrt\Lambda$ has the same adapted bases, and its constant is
$c/\Lambda<1$.  By Lemma~\ref{lem:red-KyFan}(c) for $D'$, the numbers $h,r,s,B$ formed with $D'$ satisfy
$r^2+s^2<1$, $G>0$, $\nu>0$ and $B\in K^\C_2$, so $R,S,H$ are parameters as in Definition~\ref{def:params}.
\item By the hypothesis $F_\C(B)\succeq0$, so Theorem~\ref{thm:red-square} for $D'$ gives
$\Upsilon_{D'}(\psi_1,\psi_2)\le\norm X^2$, that is, $\Upsilon(\psi_1,\psi_2)\le\Lambda\norm X^2$ by \eqref{eq:red-normal}.
\item This holds for every $\Lambda>c$, and $\Upsilon(\psi_1,\psi_2)$ does not depend on $\Lambda$;
letting $\Lambda\downarrow c$ gives $\Upsilon(\psi_1,\psi_2)\le c\norm X^2$.  Since
$X,\psi_1,\psi_2$ are arbitrary, (2V) holds.\qedhere
\end{enumerate}
\end{proof}

\begin{remark}[the boundary $G=0$, $\nu=0$]\label{rem:red-boundary}
Normalising with $\Lambda=c$ itself, that is, passing to the weight $D/\sqrt c$, can give $G=\nu=0$ (for instance if $\psi_1,\psi_2$ are eigenvectors of $D$
for $d_1,d_2$), where $F_\C$ is not defined.  The proof above uses $F_\C$ only for $\Lambda>c$ and takes the limit in
the scalar inequality $\Upsilon(\psi_1,\psi_2)\le\Lambda\norm X^2$, so $G=0$ and $\nu=0$ need no separate treatment.
\end{remark}

\section{\texorpdfstring{The case $m=2$ and the proofs of the theorems}{The case m=2 and the proofs of the theorems}}\label{sec:m2}

In this section we prove that $F_\C(B)\succeq0$ for $B\in K^\C_2$ (the case $m=2$;
Theorem~\ref{thm:m2-D}), deduce Theorems~\ref{thm:intro-norm}, \ref{cor:14} and~\ref{thm:intro-2V}, and determine the cases of
equality (Theorem~\ref{thm:intro-eq}, Section~\ref{subsec:m2-equality}).  Throughout, $m=2$, $R,S\ge0$,
$H\ge0$, $G=1-R-S-H$, and $B\in M_2(\C)$; real $B$ are the special case $F_\C(B)=F(B)$, $K_2=K^\C_2\cap M_2(\R)$,
and every formula below is then read with $\abs z^2=z^2$.  The parameters
$\nu,\sigma,\varrho,a_1,a_2,\gamma_p,\gamma_q,c_D$ and the matrix $F=F_\C(B)$ are those
of Definitions~\ref{def:params} and~\ref{def:F}.  Until Theorem~\ref{thm:m2-D}~(ii) we
only need $\nu\sigma\varrho G\ne0$; the sign of $G$ is used only in parts (ii) and (iii) of that theorem.
We write
\[
B=\begin{pmatrix}k&b_{12}\\ b_{21}&b_{22}\end{pmatrix},\qquad \det B=kb_{22}-b_{12}b_{21},
\]
\[
x_1=(S_p)_{11}=G-\abs k^2-\abs{b_{12}}^2,\qquad y_1=(S_q)_{11}=G-\abs k^2-\abs{b_{21}}^2,
\]
\[
\kappa:=G-\abs k^2=x_1+\abs{b_{12}}^2=y_1+\abs{b_{21}}^2,
\]
$d_p=\det S_p$, $d_q=\det S_q$ (real numbers, since $S_p,S_q$ are Hermitian), and, for the matrices
$\widehat S_p,\widehat S_q$ of Definition~\ref{def:F},
\[
\widehat S_p=S_p+\eta_pE_{22},\qquad
\widehat S_q=\overline{S_q}+\eta_qE_{22},
\]
where $\eta_p:=\gamma_p-a_2=S\gamma_p$ and $\eta_q:=\gamma_q-a_1=R\gamma_q$.  Here $\overline{S_q}=D_q-B^\T\bar B$ has
the same $(1,1)$ entry and the same determinant as $S_q$.  For Hermitian $2\times2$ matrices $X,Y$ put
\[
\langle X,Y\rangle_{\adj}:=\tr\bigl(\adj(X)Y\bigr)=X_{11}Y_{22}+X_{22}Y_{11}-2\operatorname{Re}(X_{12}Y_{21});
\]
this form is real, symmetric and real bilinear, $\langle X,X\rangle_{\adj}=2\det X$, and
$\det(X+Y)=\det X+\det Y+\langle X,Y\rangle_{\adj}$ (for $2\times2$ matrices $\adj$ is linear).  We use repeatedly
\[
\nu+R+S=1,\quad \nu+S=\varrho,\quad \nu+R=\sigma,\quad \sigma\varrho=\nu+RS,\quad
\sigma\gamma_p=a_2,\quad \varrho\gamma_q=a_1 .
\]

\subsection{A bordered form of \texorpdfstring{$F$}{F}}

Let $b=(se_1;re_1)$ and $\mathbf m=\mathbf m(B)=(sBe_1;rB^\T e_1)$ be as in Definition~\ref{def:F}, and put
$\mathsf U:=[\,b,\ \mathbf m\,]\in M_{4\times2}(\C)$ and $\mathsf C:=\diag(G,-1)$.  Then
$\mathsf U\mathsf C\mathsf U^*=Gbb^\T-\mathbf m\mathbf m^*$, so Definition~\ref{def:F} reads
\[
F=\widehat S_p\oplus\widehat S_q+\frac1\nu\,\mathsf U\mathsf C\mathsf U^* .
\]
In components $\mathsf U=\begin{pmatrix}sX_p\\ rX_q\end{pmatrix}$ with
$X_p=[e_1,Be_1]=\begin{pmatrix}1&k\\ 0&b_{21}\end{pmatrix}$ and
$X_q=[e_1,B^\T e_1]=\begin{pmatrix}1&k\\ 0&b_{12}\end{pmatrix}$.

\begin{lemma}[determinant lemma]\label{lem:m2-det}
Put $\mathsf A_p:=X_p^*\adj(\widehat S_p)X_p$, $\mathsf A_q:=X_q^*\adj(\widehat S_q)X_q$ (Hermitian) and
$\mathsf N:=\diag(\frac\nu G,-\nu)$.  Then
\begin{multline*}
-\frac{\nu^2}G\det F=\det\mathsf N\,\det\widehat S_p\det\widehat S_q
+S^2\abs{b_{21}}^2\det\widehat S_q+R^2\abs{b_{12}}^2\det\widehat S_p\\
+S\det\widehat S_q\,\langle\mathsf N,\mathsf A_p\rangle_{\adj}+R\det\widehat S_p\,\langle\mathsf N,\mathsf A_q\rangle_{\adj}
+RS\,\langle\mathsf A_p,\mathsf A_q\rangle_{\adj},
\end{multline*}
and $\det\mathsf N=-\frac{\nu^2}G$.
\end{lemma}

\begin{proof}
Put $\widehat S:=\widehat S_p\oplus\widehat S_q$ and
$\mathcal M:=\begin{pmatrix}\widehat S&\mathsf U\\ \mathsf U^*&-\nu\mathsf C^{-1}\end{pmatrix}$ (Hermitian).  The Schur
complement of the lower right block is $\widehat S+\frac1\nu\mathsf U\mathsf C\mathsf U^*=F$ (Definition~\ref{def:F}), so
$\det\mathcal M=\det(-\nu\mathsf C^{-1})\det F=-\frac{\nu^2}G\det F$.  If $\widehat S$ is invertible, the Schur
complement of the upper left block gives, since $\nu\mathsf C^{-1}=\mathsf N$ and $\det(-X)=\det X$ for
$2\times2$ matrices,
\begin{gather*}
\det\mathcal M=\det\widehat S\,\det(\mathsf N+\mathsf U^*\widehat S^{-1}\mathsf U),\\
\mathsf U^*\widehat S^{-1}\mathsf U=SX_p^*\widehat S_p^{-1}X_p+RX_q^*\widehat S_q^{-1}X_q
=\frac S{\det\widehat S_p}\mathsf A_p+\frac R{\det\widehat S_q}\mathsf A_q .
\end{gather*}
Expand $\det(\mathsf N+Y+Z)=\det\mathsf N+\det Y+\det Z+\langle\mathsf N,Y\rangle_{\adj}
+\langle\mathsf N,Z\rangle_{\adj}+\langle Y,Z\rangle_{\adj}$ with $Y,Z$ the two terms, use
$\det\mathsf A_p=\abs{\det X_p}^2\det\adj(\widehat S_p)=\abs{b_{21}}^2\det\widehat S_p$,
$\det\mathsf A_q=\abs{b_{12}}^2\det\widehat S_q$, and multiply by $\det\widehat S=\det\widehat S_p\det\widehat S_q$.  Both
sides are polynomials in $\operatorname{Re}B$, $\operatorname{Im}B$, and so is $\det\widehat S$, which is not identically
zero (for $B=tI_2$, $t\in\R$, it is the monic polynomial $(G-t^2)^2(\gamma_p-t^2)(\gamma_q-t^2)$); so the identity holds on
a dense set, hence everywhere.
\end{proof}

\begin{lemma}[adjugate quadratic forms]\label{lem:m2-adj}
\[
\mathsf A_p=\begin{pmatrix}\gamma_p-\abs{b_{21}}^2-\abs{b_{22}}^2&\gamma_pk-\overline{b_{22}}\det B\\
\overline{\gamma_pk-\overline{b_{22}}\det B}&\gamma_p\abs k^2+G\abs{b_{21}}^2-\abs{\det B}^2\end{pmatrix},
\]
and $\mathsf A_q$ is obtained by replacing $\gamma_p$ by $\gamma_q$ and $b_{21}$ by $b_{12}$.
\end{lemma}

\begin{proof}
For $2\times2$ matrices $\adj$ is linear, $\adj(XY)=\adj(Y)\adj(X)$ and $\adj(X^*)=\adj(X)^*$, so
$\adj(\widehat S_p)=\diag(\gamma_p,G)-\adj(B)^*\adj(B)$.  Since $\adj(B)e_1=(b_{22},-b_{21})^\T$ and
$\adj(B)Be_1=(\det B)e_1$, the Hermitian form $(\xi,\eta)\mapsto\ip{\adj(B)\xi}{\adj(B)\eta}$ takes
the values $\abs{b_{21}}^2+\abs{b_{22}}^2$ at $(e_1,e_1)$, $\overline{b_{22}}\det B$ at $(e_1,Be_1)$ and
$\abs{\det B}^2$ at $(Be_1,Be_1)$, while $\diag(\gamma_p,G)$ gives $\gamma_p$, $\gamma_pk$ and
$\gamma_p\abs k^2+G\abs{b_{21}}^2$; the columns of $X_p$ are $e_1$ and $Be_1$.  For the $q$ side,
$\widehat S_q=\diag(G,\gamma_q)-B^\T(B^\T)^*$, so $\adj(\widehat S_q)=\diag(\gamma_q,G)-\adj(B^\T)^*\adj(B^\T)$,
the columns of $X_q$ are $e_1$ and $B^\T e_1$, $\adj(B^\T)e_1=(b_{22},-b_{12})^\T$ and $\det B^\T=\det B$.
\end{proof}

\begin{lemma}[translation into state quantities]\label{lem:m2-state}
\begin{enumerate}
\item[(i)] $\det\widehat S_p=d_p+\eta_px_1$ and $\det\widehat S_q=d_q+\eta_qy_1$.
\item[(ii)] $\chi_B:=\abs{\det B}^2-G\abs{b_{22}}^2$ satisfies $\chi_B=d_p-a_2x_1+G\abs{b_{21}}^2=d_q-a_1y_1+G\abs{b_{12}}^2$.
\item[(iii)] $d_p=\det(GI-BB^*)+Rx_1$ and $d_q=\det(GI-B^*B)+Sy_1$, where
$\det(GI-BB^*)=\det(GI-B^*B)=G^2-G\norm B_2^2+\abs{\det B}^2$.  Consequently $d_p-d_q=Rx_1-Sy_1$,
\[
\Phi:=d_p+Sy_1=d_q+Rx_1=\det(GI-BB^*)+Rx_1+Sy_1,\qquad \chi_B=\Phi+G\kappa-a_1y_1-a_2x_1 .
\]
We also write $\Theta_p:=\Phi-d_p=Sy_1$ and $\Theta_q:=\Phi-d_q=Rx_1$.
\end{enumerate}
\end{lemma}

\begin{proof}
(i) For any $d$, $\det(\diag(G,d)-BB^*)=Gd-G(BB^*)_{22}-d(BB^*)_{11}+\det(BB^*)
=dx_1-G(\abs{b_{21}}^2+\abs{b_{22}}^2)+\abs{\det B}^2$.  Taking $d=\gamma_p$ and $d=a_2$ and subtracting gives the first
formula; the second is the same computation with $B^\T\bar B$ (or $B^*B$, which has the same diagonal and
determinant), for which $\det(\diag(G,d)-B^\T\bar B)=dy_1-G(\abs{b_{12}}^2+\abs{b_{22}}^2)+\abs{\det B}^2$.
(ii) Solve the two expressions of (i) with $d=a_2$, $d=a_1$ for $\abs{\det B}^2-G\abs{b_{22}}^2$.
(iii) $D_p=GI+RE_{22}$, and $\det(X+tE_{22})=\det X+tX_{11}$ for $2\times2$ matrices $X$, with
$(GI-BB^*)_{11}=x_1$; likewise $D_q=GI+SE_{22}$ and $(GI-B^*B)_{11}=y_1$.  The two determinants agree because
$BB^*$ and $B^*B$ have the same eigenvalues, and the closed form is the formula in the proof of (i) with
$d=G$.  Finally $\chi_B=d_p-a_2x_1+G(\kappa-y_1)=\Phi-Sy_1-a_2x_1+G\kappa-Gy_1$.
\end{proof}

\begin{lemma}[linear border terms and the cross term]\label{lem:m2-cross}
Put $P:=\abs{b_{12}}^2$ and $Q:=\abs{b_{21}}^2$.
\begin{enumerate}
\item[(i)] $SQ+\langle\mathsf N,\mathsf A_p\rangle_{\adj}=(\nu+S)Q-\frac\nu G\bigl(\gamma_pP+\det\widehat S_p\bigr)$ and
$RP+\langle\mathsf N,\mathsf A_q\rangle_{\adj}=(\nu+R)P-\frac\nu G\bigl(\gamma_qQ+\det\widehat S_q\bigr)$.
\item[(ii)] $\langle\mathsf A_p,\mathsf A_q\rangle_{\adj}=\gamma_pPy_1+\gamma_qQx_1-2GPQ+(P+Q)\chi_B$.
\end{enumerate}
\end{lemma}

\begin{proof}
(i) Since $\adj\mathsf N=\diag(-\nu,\frac\nu G)$, $\langle\mathsf N,\mathsf A_p\rangle_{\adj}=-\nu(\mathsf A_p)_{11}+\frac\nu G(\mathsf A_p)_{22}$;
only the diagonal entries of $\mathsf A_p$ enter.
By Lemma~\ref{lem:m2-state}~(i) (with the formula in its proof),
$\abs{\det B}^2=\det\widehat S_p-\gamma_px_1+G(\abs{b_{21}}^2+\abs{b_{22}}^2)$, and $\abs k^2+x_1=G-\abs{b_{12}}^2$; hence
$(\mathsf A_p)_{22}=\gamma_pG-\gamma_pP-\det\widehat S_p-G\abs{b_{22}}^2$.  With
$(\mathsf A_p)_{11}=\gamma_p-Q-\abs{b_{22}}^2$, the terms $\nu\gamma_p$ and $\nu\abs{b_{22}}^2$ cancel and
$\langle\mathsf N,\mathsf A_p\rangle_{\adj}=\nu Q-\frac\nu G(\gamma_pP+\det\widehat S_p)$; add $SQ$.  The $q$ side is identical.

(ii) Write also $W:=\abs{b_{22}}^2$ and $\mathsf D:=\abs{\det B}^2$, so that
$\chi_B=\mathsf D-GW$.  By Lemma~\ref{lem:m2-adj},
\begin{multline*}
\langle\mathsf A_p,\mathsf A_q\rangle_{\adj}=(\gamma_p-Q-W)(\gamma_q\abs k^2+GP-\mathsf D)+(\gamma_p\abs k^2+GQ-\mathsf D)(\gamma_q-P-W)\\
-2\operatorname{Re}\bigl((\gamma_pk-\overline{b_{22}}\det B)\,\overline{(\gamma_qk-\overline{b_{22}}\det B)}\bigr).
\end{multline*}
This is the only place where the phases of the entries of $B$ enter, and the last term equals
$-2\gamma_p\gamma_q\abs k^2+2(\gamma_p+\gamma_q)\operatorname{Re}\bigl(kb_{22}\overline{\det B}\bigr)-2W\mathsf D$.
Taking the squared modulus of $b_{12}b_{21}=kb_{22}-\det B$ gives
\[
2\operatorname{Re}\bigl(kb_{22}\overline{\det B}\bigr)=\mathsf D+\abs k^2W-PQ ,
\]
so the phase is absorbed into $\mathsf D=\abs{\det B}^2$, and the last term equals
$-2\gamma_p\gamma_q\abs k^2+(\gamma_p+\gamma_q)(\mathsf D+\abs k^2W-PQ)-2W\mathsf D$.  Collecting terms: the multiples of
$\gamma_p\gamma_q\abs k^2$ ($1+1-2$) and of $W\mathsf D$ ($1+1-2$) cancel; the terms with $\gamma_p$ are
$\gamma_pGP-\gamma_p\abs k^2P-\gamma_pPQ=\gamma_pP\,(G-\abs k^2-Q)=\gamma_pPy_1$; those with $\gamma_q$ are
$\gamma_qQx_1$ in the same way; the rest is $-2GPQ+(P+Q)(\mathsf D-GW)$.  This is (ii).
\end{proof}

\subsection{The four-quantity identity}

\begin{proposition}[four-quantity identity]\label{prop:m2-four}
Assume $\nu\sigma\varrho G\ne0$.  For every $B\in M_2(\C)$,
\begin{multline*}
\nu^2\sigma\varrho\det F=\sigma\varrho\nu\Phi^2-H\sigma\varrho(Rx_1+Sy_1)\Phi+H^2RSx_1y_1\\
+HRS\kappa\bigl[(2\nu+R+S)\Phi-H(Rx_1+Sy_1)\bigr].
\end{multline*}
\end{proposition}

\begin{proof}
Write $\widehat P_p:=\det\widehat S_p$ and $\widehat P_q:=\det\widehat S_q$.  In Lemma~\ref{lem:m2-det}
with $P=\abs{b_{12}}^2$ and $Q=\abs{b_{21}}^2$ as in Lemma~\ref{lem:m2-cross}, group the second and fourth terms as
$S\widehat P_q\bigl(SQ+\langle\mathsf N,\mathsf A_p\rangle_{\adj}\bigr)$
and the third and fifth as $R\widehat P_p\bigl(RP+\langle\mathsf N,\mathsf A_q\rangle_{\adj}\bigr)$, and
insert Lemma~\ref{lem:m2-cross}~(i).  The coefficient of $\widehat P_p\widehat P_q$ becomes
$-\frac{\nu^2}G-\frac\nu G(R+S)=-\frac\nu G$, so
\[
-\frac{\nu^2}G\det F=-\frac\nu G\widehat P_p\widehat P_q
+S\widehat P_q\Bigl[(\nu+S)Q-\frac\nu G\gamma_pP\Bigr]
+R\widehat P_p\Bigl[(\nu+R)P-\frac\nu G\gamma_qQ\Bigr]+RS\langle\mathsf A_p,\mathsf A_q\rangle_{\adj} .
\]
Multiply by $-G\sigma\varrho$ and use $\nu+S=\varrho$, $\nu+R=\sigma$, $\sigma\gamma_p=a_2$,
$\varrho\gamma_q=a_1$ and Lemma~\ref{lem:m2-cross}~(ii).  Then $\nu^2\sigma\varrho\det F$ is the sum of
\begin{align*}
\text{part 1}&=\nu\,(\sigma\widehat P_p)(\varrho\widehat P_q),\\
\text{part 2}&=-S\,(\varrho\widehat P_q)\bigl(G\sigma\varrho Q-\nu a_2P\bigr),\\
\text{part 3}&=-R\,(\sigma\widehat P_p)\bigl(G\sigma\varrho P-\nu a_1Q\bigr),\\
\text{part 4}&=-RS\Bigl[G\varrho a_2Py_1+G\sigma a_1Qx_1
+G\sigma\varrho\bigl(-2GPQ+(P+Q)\chi_B\bigr)\Bigr].
\end{align*}
By Lemma~\ref{lem:m2-state} (with $\eta_p=S\gamma_p$, $\eta_q=R\gamma_q$),
\[
\sigma\widehat P_p=\sigma\Phi-\sigma Sy_1+Sa_2x_1,\quad \varrho\widehat P_q=\varrho\Phi-\varrho Rx_1+Ra_1y_1,\quad
P=\kappa-x_1,\quad Q=\kappa-y_1,
\]
and $\chi_B=\Phi+G\kappa-a_1y_1-a_2x_1$; these are linear forms in $(\Phi,\kappa,x_1,y_1)$, so each part is a
quadratic form in these four quantities.  The exchange of $p$ and $q$ ($R\leftrightarrow S$, $x_1\leftrightarrow y_1$,
$P\leftrightarrow Q$, hence $\sigma\leftrightarrow\varrho$, $a_1\leftrightarrow a_2$) fixes parts 1 and 4 and $\chi_B$ and
interchanges parts 2 and 3; so the rows $\Phi y_1$, $\kappa y_1$, $y_1^2$ are the mirror images of the rows $\Phi x_1$,
$\kappa x_1$, $x_1^2$.  The monomial $\Phi^2$ occurs only in part 1, with coefficient $\sigma\varrho\nu$, and $\kappa^2$
does not occur: only in part 4 do both factors contain $\kappa$, and there $-2GPQ$ and $(P+Q)\chi_B$ give
$-2G\kappa^2$ and $+2G\kappa^2$.  The other coefficients, obtained by multiplying the coefficients of the two linear
factors, are listed below.  For example, $(\Phi,x_1)$ has the coefficients $(\sigma,Sa_2)$ in $\sigma\widehat P_p$ and
$(\varrho,-\varrho R)$ in $\varrho\widehat P_q$, so the entry of part~1 in the row $\Phi x_1$ is
$\nu\bigl(\sigma\cdot(-\varrho R)+Sa_2\cdot\varrho\bigr)=\nu\varrho(Sa_2-\sigma R)$.  In parts 2 and 3 the coefficients of $\kappa$ in the second factor are
$G\sigma\varrho-\nu a_2=-R(\nu\sigma+HS)$ and $G\sigma\varrho-\nu a_1=-S(\nu\varrho+HR)$, by $a_2=\sigma-H$, $a_1=\varrho-H$,
$G=\nu-H$, $1-\varrho=R$, $1-\sigma=S$ and $\sigma\varrho-\nu=RS$.

\providecommand{\mtwotablesize}{\footnotesize}
\medskip
\begin{center}\mtwotablesize\setlength{\tabcolsep}{3pt}
\begin{tabular}{@{}llllll@{}}
\toprule
 & part 1 & part 2 & part 3 & part 4 & sum\\
\midrule
$\Phi\kappa$ & $0$ & $RS\varrho(\nu\sigma+HS)$ & $RS\sigma(\nu\varrho+HR)$ & $-2RSG\sigma\varrho$ & $HRS(2\nu+R+S)$\\
$\Phi x_1$ & $\nu\varrho(Sa_2-\sigma R)$ & $-\nu S\varrho a_2$ & $R\sigma^2\varrho G$ & $RS\sigma\varrho G$ & $-HR\sigma\varrho$\\
$\kappa x_1$ & $0$ & $-R^2S\varrho(\nu\sigma+HS)$ & $RS^2a_2(\nu\varrho+HR)$ & $RSG\sigma(\varrho(R-S)+HR)$ & $-H^2R^2S$\\
$x_1^2$ & $-\nu RS\varrho a_2$ & $\nu RS\varrho a_2$ & $RS\sigma\varrho Ga_2$ & $-RS\sigma\varrho Ga_2$ & $0$\\
$x_1y_1$ & $\nu RS(\sigma\varrho+a_1a_2)$ & $-RS(G\sigma\varrho^2+\nu a_1a_2)$ & $-RS(G\sigma^2\varrho+\nu a_1a_2)$
 & $(\ast)$ & $H^2RS$\\
\bottomrule
\end{tabular}
\end{center}
\noindent where $(\ast)=RSG\bigl[\varrho a_2+\sigma a_1+2G\sigma\varrho-\sigma\varrho(a_1+a_2)\bigr]$.
\medskip

\noindent The sums are checked as follows; the $x_1^2$ entries cancel in pairs.
\begin{itemize}
\item $\Phi\kappa$: by $\nu-G=H$ the row is $HRS(2\sigma\varrho+S\varrho+R\sigma)$, and
$2\sigma\varrho+S\varrho+R\sigma=2(\nu+RS)+(S-RS)+(R-RS)=2\nu+R+S$.
\item $\Phi x_1$: parts 1 and 2 give $-\nu\sigma\varrho R$, parts 3 and 4 give $R\sigma\varrho G(\sigma+S)=R\sigma\varrho G$,
and $G-\nu=-H$.
\item $\kappa x_1$: divide by $RS$ and insert $a_2=\sigma-H$, $G=\nu-H$.  Sorted by powers of $H$, parts 2, 3, 4 give
$-R\nu\sigma\varrho-HRS\varrho$, $\;S\nu\sigma\varrho+H(RS\sigma-S\nu\varrho)-H^2RS$ and
$(R-S)\nu\sigma\varrho+H(R\nu\sigma-(R-S)\sigma\varrho)-H^2R\sigma$.  The $H^0$ terms cancel; the $H^1$ terms add to
$RS(\sigma-\varrho)+\nu(R\sigma-S\varrho)-(R-S)\sigma\varrho=(R-S)(RS+\nu-\sigma\varrho)=0$, since
$\sigma-\varrho=R\sigma-S\varrho=R-S$; the $H^2$ terms give $-R(S+\sigma)=-R$.
\item $x_1y_1$: divide by $RS$ and put $\mathsf p:=\sigma\varrho$, $\mathsf s:=\sigma+\varrho=1+\nu$, so that
$a_1a_2=\mathsf p-H\mathsf s+H^2$, $\varrho a_2+\sigma a_1=2\mathsf p-H\mathsf s$ and $a_1+a_2=\mathsf s-2H$.  The sum of the
four entries is then
$\nu\mathsf p-\nu a_1a_2-2G\mathsf p\mathsf s+2G\mathsf p-GH\mathsf s+2G^2\mathsf p+2GH\mathsf p$, and with
$\nu\mathsf p-\nu a_1a_2=\nu H\mathsf s-\nu H^2$, $2G^2\mathsf p+2GH\mathsf p=2G\nu\mathsf p$ and $\nu H\mathsf s-GH\mathsf s=H^2\mathsf s$
it equals $H^2(\mathsf s-\nu)+2G\mathsf p(1+\nu-\mathsf s)=H^2$.
\end{itemize}
With the mirror rows ($-HS\sigma\varrho$, $-H^2RS^2$, $0$), these are the coefficients of the right side of the
proposition.
\end{proof}

\subsection{The determinant identity and positivity}

\begin{theorem}\label{thm:m2-D}
Let $H\ge0$ and put $\Pi_p:=\frac{G\Phi+Hd_p}\nu$ and $\Pi_q:=\frac{G\Phi+Hd_q}\nu$.
\begin{enumerate}
\item[(i)] \textup{(Identity.)}  If $\nu\sigma\varrho G\ne0$, then for every $B\in M_2(\C)$
\[
\nu^3\sigma\varrho\det F=\sigma\varrho\,(G\Phi+Hd_p)(G\Phi+Hd_q)+HRS\,\mathcal R_2,
\]
where
\begin{align*}
\mathcal R_2&=\nu\kappa(2\nu+R+S)\Phi-H\nu\kappa(\Theta_p+\Theta_q)-H\Theta_p\Theta_q\\
&=\kappa\Phi\,c_D+H\nu\kappa(d_p+d_q)+H\Theta_pd_q+HS\abs{b_{21}}^2\Phi .
\end{align*}
\item[(ii)] If $G>0$ and $B\in K^\C_2$, then $\nu\det F\ge\Pi_p\Pi_q\ge d_pd_q\ge0$.
\item[(iii)] If $G>0$ and $B\in K^\C_2$, then $F_\C(B)\succeq0$.  In particular $F(B)\succeq0$ for $B\in K_2$.
\end{enumerate}
\end{theorem}

\begin{proof}
(i) By Lemma~\ref{lem:m2-state}, $Rx_1=\Theta_q$, $Sy_1=\Theta_p$ and
$H^2RSx_1y_1=H^2\Theta_p\Theta_q$, so Proposition~\ref{prop:m2-four} reads
\[
\nu^2\sigma\varrho\det F=\sigma\varrho\nu\Phi^2-H\sigma\varrho(\Theta_p+\Theta_q)\Phi+H^2\Theta_p\Theta_q
+HRS\kappa\bigl[(2\nu+R+S)\Phi-H(\Theta_p+\Theta_q)\bigr].
\]
Complete the product:
$\sigma\varrho\nu\Phi^2-H\sigma\varrho(\Theta_p+\Theta_q)\Phi
=\sigma\varrho\nu\bigl(\Phi-\frac H\nu\Theta_p\bigr)\bigl(\Phi-\frac H\nu\Theta_q\bigr)-\sigma\varrho\frac{H^2}\nu\Theta_p\Theta_q$,
and $\Phi-\frac H\nu\Theta_p=\frac{(\nu-H)\Phi+H(\Phi-\Theta_p)}\nu=\frac{G\Phi+Hd_p}\nu=\Pi_p$; likewise for $q$.
The total coefficient of $\Theta_p\Theta_q$ is $H^2-\sigma\varrho\frac{H^2}\nu=-\frac{H^2RS}\nu$ because
$\sigma\varrho=\nu+RS$.  Hence
\[
\nu^2\sigma\varrho\det F=\sigma\varrho\nu\,\Pi_p\Pi_q+\frac{HRS}\nu\Bigl[\nu\kappa(2\nu+R+S)\Phi-H\nu\kappa(\Theta_p+\Theta_q)-H\Theta_p\Theta_q\Bigr],
\]
and multiplying by $\nu$ gives the identity with the first expression of $\mathcal R_2$.  For the second
expression use $\nu(2\nu+R+S)=c_D+2H\nu+HS$, $S\kappa-\Theta_p=S(\kappa-y_1)=S\abs{b_{21}}^2$,
$\Phi-\Theta_p=d_p$ and $\Phi-\Theta_q=d_q$: then
\begin{align*}
&\kappa\Phi c_D+H\nu\kappa(\Phi-\Theta_p)+H\nu\kappa(\Phi-\Theta_q)+H\Theta_p(\Phi-\Theta_q)+H\Phi(S\kappa-\Theta_p)\\
&\qquad=\kappa\Phi(c_D+2H\nu+HS)-H\nu\kappa(\Theta_p+\Theta_q)-H\Theta_p\Theta_q ,
\end{align*}
which is the first expression.

(ii) On $K^\C_2$ we have $d_p,d_q\ge0$ (determinants of positive semidefinite matrices),
$x_1,y_1\ge0$ (diagonal entries), hence $\Theta_p=Sy_1\ge0$, $\Theta_q=Rx_1\ge0$,
$\kappa=x_1+\abs{b_{12}}^2\ge0$, $\Phi=d_p+\Theta_p\ge0$, and $c_D=2\nu G+\nu R+SG>0$ does not
depend on $B$.  So every term of the second expression of $\mathcal R_2$ is a product of nonnegative factors
and $\mathcal R_2\ge0$.  Since $\nu,\sigma,\varrho>0$, dividing (i) by $\nu^3\sigma\varrho$ gives
$\nu\det F\ge\Pi_p\Pi_q$.  Finally $\Pi_p=\frac G\nu\Phi+\frac H\nu d_p$ is a convex combination ($G+H=\nu$) of
$\Phi=d_p+\Theta_p\ge d_p$ and $d_p$, so $\Pi_p\ge d_p\ge0$, and likewise $\Pi_q\ge d_q\ge0$.

(iii) Assume first that $S_p\succ0$ and $S_q\succ0$.  Then $\widehat S_p=S_p+\eta_pE_{22}\succ0$ and
$\widehat S_q=\overline{S_q}+\eta_qE_{22}\succ0$, since $\eta_p,\eta_q\ge0$.  Since
$\mathsf U\mathsf C\mathsf U^*=Gbb^\T-\mathbf m\mathbf m^*$, Definition~\ref{def:F} reads
\[
F=\widehat S_p\oplus\widehat S_q+\frac G\nu\,bb^\T-\frac1\nu\,\mathbf m\mathbf m^* .
\]
For $\mathbf x\ne0$ with $\mathbf m^*\mathbf x=0$ we get
$\mathbf x^*F\mathbf x\ge\mathbf x^*(\widehat S_p\oplus\widehat S_q)\mathbf x>0$, because $G>0$.  These $\mathbf x$,
together with $0$, form a subspace of dimension at least $3$, so by the Courant--Fischer minimax principle $F$ has at
least three positive eigenvalues.  By (ii), $\det F\ge\nu^{-1}d_pd_q>0$, and $\det F$ is the product of the four
eigenvalues; so the fourth one is positive as well, and $F\succ0$.  For general $B\in K^\C_2$ and $\tau\in[0,1)$,
\[
S_p(\tau B)=D_p-\tau^2BB^*=(1-\tau^2)D_p+\tau^2S_p(B)\succeq(1-\tau^2)D_p\succ0,
\]
because $D_p=\diag(G,a_2)\succ0$, and likewise $S_q(\tau B)\succ0$.  Hence $F_\C(\tau B)\succ0$, and letting
$\tau\uparrow1$ gives $F_\C(B)\succeq0$, by continuity.  For real $B\in K_2$, $F_\C(B)=F(B)$ and $B\in K^\C_2$.
\end{proof}

\begin{remark}\label{rem:m2-remarks}
The cases $H=0$ (where $\nu\det F=\Phi^2$), $R=0$ and $S=0$ are contained in Theorem~\ref{thm:m2-D}
directly, and (iii) holds on the whole closed set $K^\C_2$, also at points with $d_p=0$, which do not arise in
Theorem~\ref{thm:red-chain} but may occur in Section~\ref{subsec:m2-equality}, where $\Lambda=c$.
\end{remark}

\subsection{\texorpdfstring{Proofs of Theorems~\ref{thm:intro-norm}, \ref{thm:intro-2V} and~\ref{cor:14}}{Proofs of the main results and of the second theorem}}\label{subsec:m2-proof14}

\begin{corollary}[all $Z$]\label{cor:m2-2V-all}
Let $n\ge2$ and $Z\in M_n(\C)$.  For all $X,Y\in M_n(\C)$ with $Y\ne0$, with $\mathsf W_Z=Z^*X^*Y+YX^*Z^*$,
\[
\norm{XZY-YZX}^2+\norm{\mathsf W_Z}^2-\frac{\abs{\ip Y{\mathsf W_Z}}^2}{\norm Y^2}\le\norm Z_{(2),2}^2\norm X^2\norm Y^2 ,
\]
and the constant $\norm Z_{(2),2}^2$ is best possible.  For $Z=\omega^{-1/2}$ this is Theorem~\ref{thm:intro-2V}.
\end{corollary}

\begin{proof}[Proof of Theorems~\ref{thm:intro-norm}, \ref{thm:intro-2V} and Corollary~\ref{cor:m2-2V-all}]
\emph{Square matrices.}  Let $n\ge2$, $a\ge b>0$, $D=\diag(a,b,\dots,b)$ and $c=a^2+b^2=\norm D_{(2),2}^2$.  By
Theorem~\ref{thm:m2-D}~(iii) the hypothesis of Theorem~\ref{thm:red-chain} holds, so (2V) holds for every $3\times3$
weight, in particular for $\diag(a,b,b)$, and by Corollary~\ref{cor:red-local} for $D$.  Proposition~\ref{prop:red-2V}
gives $\mathcal L_D(X,Y)\le c\norm X^2\norm Y^2$ for all $X,Y$ with $Y\ne0$, that is, $\phi(D)\le\norm D_{(2),2}$, and
hence also $\mathfrak n(D)\le\norm D_{(2),2}$, since $\mathfrak n\le\phi$.  Proposition~\ref{prop:red-two-level}, with
$\theta=\phi$ and with $\theta=\mathfrak n$, gives $\phi(Z)=\mathfrak n(Z)=\norm Z_{(2),2}$ for every $Z\in M_n(\C)$.  The
left side of the inequality of Corollary~\ref{cor:m2-2V-all} is $\mathcal L_Z(X,Y)$, so this is
Corollary~\ref{cor:m2-2V-all}, and it is \eqref{eq:intro-norm} for $p=q=n$ together with the last statement of
Theorem~\ref{thm:intro-norm} (for $n=1$ the left side of \eqref{eq:intro-norm} is $0$).  For
Theorem~\ref{thm:intro-2V} take $Z=\omega^{-1/2}$: $Z$ is Hermitian, so $\mathsf W_Z=\mathsf W$, and
$\norm Z_{(2),2}^2=\lambda_1^{-1}+\lambda_2^{-1}$.

\emph{Rectangular matrices.}  Let $N=\max(p,q)$, and for $k\le N$ let $\iota_k\in M_{N\times k}(\C)$ be the
inclusion of $\C^k$ as the first $k$ coordinates, so that $\iota_k^*\iota_k=I_k$.  Put $X_N=\iota_pX\iota_q^*$,
$Y_N=\iota_pY\iota_q^*$ and $Z_N=\iota_qZ\iota_p^*$ in $M_N(\C)$.  Then
$X_NZ_NY_N-Y_NZ_NX_N=\iota_p(XZY-YZX)\iota_q^*$, the Frobenius norms are unchanged,
and $Z_N$ has the nonzero singular values of $Z$, so $\norm{Z_N}_{(2),2}=\norm Z_{(2),2}$; the square case
gives \eqref{eq:intro-norm}.

\emph{Equality.}  Let $p,q\ge2$, and let $Z=U\Sigma V$ with unitaries $U\in M_q(\C)$, $V\in M_p(\C)$ and
$\Sigma\in M_{q\times p}(\R)$ having the entries $\Sigma_{kk}=s_k(Z)$ and zeros elsewhere.  Let $E_{12},E_{21}\in M_{p\times q}(\C)$
be matrix units, and put $X=V^*E_{12}U^*$, $Y=V^*E_{21}U^*$.  Then $XZY=V^*E_{12}\Sigma E_{21}U^*=s_2(Z)V^*E_{11}U^*$ and
$YZX=V^*E_{21}\Sigma E_{12}U^*=s_1(Z)V^*E_{22}U^*$, so $\norm{XZY-YZX}^2=s_1(Z)^2+s_2(Z)^2=\norm Z_{(2),2}^2\norm X^2\norm Y^2$.
Let $p=1<q$, so that $Z=z\in\C^q$; for $z=0$ every pair gives equality.  For $z\ne0$ take $X=z^*/\norm z$ and
$Y=w^*$ with a unit vector $w\perp z$.  Then $XZ=\norm z$ and $YZ=0$, so $XZY-YZX=\norm z\,Y$ and
$\norm{XZY-YZX}=\norm z=\norm Z_{(2),2}\norm X\norm Y$.  The case $q=1<p$ follows by transposition.
\end{proof}

\begin{proof}[Proof of Theorem~\ref{cor:14}]
We prove part (ii); part (i) is its special case of real $A,B$ and real symmetric $\omega$.  With $X=A\omega^{1/2}$,
$Y=B\omega^{1/2}$ and $Z=\omega^{-1/2}$ one has $\norm X=\norm A_\omega$, $\norm Y=\norm B_\omega$ and
$XZY-YZX=[A,B]\,\omega^{1/2}$ (Lemma~\ref{lem:red-normal}), and $\norm Z_{(2),2}^2=\lambda_1^{-1}+\lambda_2^{-1}$, so
\eqref{eq:intro-norm} is \eqref{eq:14}.  The constant is best possible by the example after
Conjecture~\ref{conj:MKOC}.
\end{proof}

\subsection{\texorpdfstring{Equality in \eqref{eq:14}}{Equality in (14)}}\label{subsec:m2-equality}

We prove Theorem~\ref{thm:intro-eq}.  For a subspace $P$ let $\Pi_P$ be the orthogonal projection onto $P$ and
$M(P)=\{Z\in M_n(\C):Z=\Pi_PZ\Pi_P\}$.  With $D=\diag(d_1\ge\dots\ge d_n>0)$ and $c=d_1^2+d_2^2$ as in
Section~\ref{sec:reduction}, we call $(X,Y)$ an \emph{equality pair for $D$} if $X,Y\ne0$ and
$\norm{T_D(X,Y)}^2=c\norm X^2\norm Y^2$.  By the proof of Lemma~\ref{lem:red-normal}, if $\omega$ is diagonal,
then $(A,B)$ gives equality in \eqref{eq:14} if and only if $(A\omega^{1/2},B\omega^{1/2})$ is an equality pair
for $D=\omega^{-1/2}$; and for a subspace $P$ invariant under $\omega$, $A\in M(P)$ if and only if
$A\omega^{1/2}\in M(P)$.

\begin{proposition}[the case $n=2$]\label{prop:eq-n2}
Let $\omega\in M_2(\C)$ be positive definite, $\varphi=\Tr(\omega\,\cdot\,)$ and $A_0=A-\varphi(A)\varphi(1)^{-1}1$,
$B_0=B-\varphi(B)\varphi(1)^{-1}1$.  Then
\[
\norm{[A,B]}_\omega^2=\Tr(\omega^{-1})\bigl(\norm{A_0}_\omega^2\norm{B_0}_\omega^2-\abs{\varphi(A_0^*B_0)}^2\bigr).
\]
Consequently \eqref{eq:14} holds, and for $A,B\ne0$ it is an equality if and only if
$\varphi(A)=\varphi(B)=\varphi(A^*B)=0$.
\end{proposition}

\begin{proof}
We may take $\omega=\diag(\lambda_1,\lambda_2)$; then $\Tr(\omega^{-1})=c$.  The matrices $u_1=E_{12}/\sqrt{\lambda_2}$,
$u_2=E_{21}/\sqrt{\lambda_1}$ and $u_3=(\lambda_2E_{11}-\lambda_1E_{22})/\sqrt{\lambda_1\lambda_2(\lambda_1+\lambda_2)}$ form an
orthonormal basis of $\{Z:\varphi(Z)=0\}$ for $\ip XY_\varphi=\varphi(X^*Y)$.  One computes
$[u_1,u_2]=(E_{11}-E_{22})/\sqrt{\lambda_1\lambda_2}$, $[u_3,u_1]=\sqrt{c}\,u_1$ and $[u_3,u_2]=-\sqrt c\,u_2$; these three
matrices have $\norm\cdot_\omega^2=c$ and are mutually orthogonal (one is diagonal, one is a multiple of $E_{12}$,
one of $E_{21}$).  Write $A_0=\sum_i\alpha_iu_i$ and $B_0=\sum_i\beta_iu_i$.  Since $[A,B]=[A_0,B_0]$,
\[
\norm{[A,B]}_\omega^2=c\sum_{i<j}\abs{\alpha_i\beta_j-\alpha_j\beta_i}^2
=c\Bigl(\norm\alpha^2\norm\beta^2-\Bigl\lvert\sum_i\bar\alpha_i\beta_i\Bigr\rvert^2\Bigr)
\]
by the Lagrange identity, which is the claim.  Since $\norm A_\omega^2=\norm{A_0}_\omega^2+\abs{\varphi(A)}^2/\varphi(1)$,
the right side is at most $c\norm A_\omega^2\norm B_\omega^2$.  For $A,B\ne0$, equality forces
$\norm{A_0}_\omega\norm{B_0}_\omega=\norm A_\omega\norm B_\omega>0$, hence $\varphi(A)=\varphi(B)=0$, and then
$\varphi(A_0^*B_0)=\varphi(A^*B)=0$; the converse is read off from the identity.
\end{proof}

The sufficiency in Theorem~\ref{thm:intro-eq} follows for every $n$: if $P$, $A$, $B$ are as there, then $[A,B]\in M(P)$,
all norms are computed on $P$ with $\omega|_P$, and $\Tr((\omega|_P)^{-1})=\lambda_1^{-1}+\lambda_2^{-1}$.

\begin{lemma}[leakage]\label{lem:eq-leak}
Let $P$ be a subspace invariant under $D$ and let $(X,Y)$ be an equality pair for $D$ with $X\in M(P)$.  Then $Y\in M(P)$.
\end{lemma}

\begin{proof}
Put $\Pi=\Pi_P$, $\Pi^\perp=1-\Pi$, $Y_{PP}=\Pi Y\Pi$, $Y_{PQ}=\Pi Y\Pi^\perp$, $Y_{QP}=\Pi^\perp Y\Pi$ and
$Y_{QQ}=\Pi^\perp Y\Pi^\perp$.  Since $X=\Pi X\Pi$ and $\Pi D=D\Pi$, we have $XD\Pi^\perp=0=\Pi^\perp DX$, so
$T_D(X,Y)=T_D(X,Y_{PP})+XDY_{PQ}-Y_{QP}DX$, and the three terms lie in $\Pi M_n\Pi$, $\Pi M_n\Pi^\perp$ and
$\Pi^\perp M_n\Pi$, which are mutually orthogonal.  By Theorem~\ref{thm:intro-norm} with $Z=D$ (inequality~\eqref{eq:red-TD} for $D$, as $\norm D_{(2),2}^2=c$) and
$\norm{XD}_{\mathrm{op}}\le d_1\norm X$,
\[
\norm{T_D(X,Y)}^2\le c\norm X^2\norm{Y_{PP}}^2+d_1^2\norm X^2\bigl(\norm{Y_{PQ}}^2+\norm{Y_{QP}}^2\bigr),
\]
so $c\norm X^2\norm Y^2-\norm{T_D(X,Y)}^2\ge c\norm X^2\norm{Y_{QQ}}^2+d_2^2\norm X^2(\norm{Y_{PQ}}^2+\norm{Y_{QP}}^2)$.
The left side is $0$, and $X\ne0$, so $Y_{PQ}=Y_{QP}=Y_{QQ}=0$.
\end{proof}

In the next lemma and proposition, $n\ge3$, $D=\diag(a,b,\dots,b)$ with $a>b>0$, $e=\mathbf e_1$, $c=a^2+b^2$, and
$\Upsilon$ is the functional of Section~\ref{sec:reduction} for a fixed $X$.  By the proof of
Theorem~\ref{thm:intro-norm} (square matrices), (2V) holds for this $D$ and for every $3\times3$ weight.
Recall $D\psi=b\psi+(a-b)\ip e\psi e$, so that
$\norm{D\psi}^2=b^2+(a^2-b^2)\abs{\ip e\psi}^2$ and, for a plane $W$, $\Tr(\Pi_WD^2)=2b^2+(a^2-b^2)\norm{\Pi_We}^2\le c$,
with equality if and only if $e\in W$.

\begin{lemma}[equality in (2V)]\label{lem:eq-2V}
Let $X\ne0$ and let $\psi_1,\psi_2$ be unit vectors with $\Upsilon(\psi_1,\psi_2)=c\norm X^2$.  Then either
\textup{(A)} $P:=\operatorname{span}\{e,\psi_1,\psi_2\}$ is two-dimensional and $X\in M(P)$, or \textup{(B)}
$\psi_1,\psi_2\in\C e$ and $X=ew^*+w'e^*$ with $w,w'\perp e$.
\end{lemma}

\begin{proof}
\emph{Case $\psi_1,\psi_2\in\C e$.}  Since $\Upsilon(\psi_1,\psi_2)$ depends only on $\psi_1\psi_1^*$ and
$\psi_2\psi_2^*$, we may take $\psi_1=\psi_2=e$.  Write $x_c=\Pi_{e^\perp}Xe$ and
$x_r=\Pi_{e^\perp}X^*e$.  Since $De=ae$ and $D=b$ on $e^\perp$, Lemma~\ref{lem:red-decomp} gives
$T_D(X,ee^*)=a(x_ce^*-ex_r^*)$ and $\Upsilon(e,e)=a^2(\norm{x_c}^2+\norm{x_r}^2)+b^2\norm{x_r}^2+b^2\norm{x_c}^2
=c(\norm{x_c}^2+\norm{x_r}^2)$.  This equals $c\norm X^2$ only if $\ip e{Xe}=0$ and $\Pi_{e^\perp}X\Pi_{e^\perp}=0$, which is (B).

\emph{Localization.}  From now on $\psi_1,\psi_2$ are not both in $\C e$.  Let $P=\operatorname{span}\{e,\psi_1,\psi_2\}$,
of dimension $2$ or $3$, and use the notation of Lemma~\ref{lem:red-local} with $V=P$.  By (2V) for $D$,
$\Upsilon_{X_{VV}}(\psi_1,\psi_2)\le c\norm{X_{VV}}^2$, so Lemma~\ref{lem:red-local} gives
\[
c\norm X^2=\Upsilon_X(\psi_1,\psi_2)\le c\norm{X_{VV}}^2+\Xi_{VW}+\Xi_{WV}
\le c\bigl(\norm X^2-\norm{X_{WW}}^2\bigr).
\]
Hence $X_{WW}=0$, $\Xi_{VW}=c\norm{X_{VW}}^2$ and $\Xi_{WV}=c\norm{X_{WV}}^2$, and by Lemma~\ref{lem:red-local}~(ii)
$X_{VW}=X_{WV}=0$; that is, $X\in M(P)$.  If $\dim P=2$, this is (A).  It remains to show that $\dim P=3$ is
impossible.

\emph{The case $\dim P=3$.}  Suppose $\dim P=3$.  As in the proof of Corollary~\ref{cor:red-local}, a unitary $U$
with $U^*DU=D$ and $U^*P=\C^3\oplus0$ moves $X$, $\psi_1$, $\psi_2$ into $\C^3\oplus0$ without changing $\Upsilon$,
$\norm X$ and $c$; so we may assume $n=3$, $D=\diag(a,b,b)$ and $P=\C^3$, and (2V) holds for this $D$.  We use the
adapted bases and the notation of Sections~\ref{subsec:red-adapted}--\ref{subsec:red-assembly}, and put
$\nu_c=c-r^2-s^2\ge0$ (Lemma~\ref{lem:red-KyFan}(b)).

If $\nu_c=0$, then $h=0$ and $r^2+s^2=c$.  With $\alpha_k=\ip e{\psi_k}$, $\norm{D\psi_2}^2=r^2$ and
$\norm{D\psi_1}^2=s^2$ (Lemma~\ref{lem:red-zero} and its proof) give $\abs{\alpha_1}^2+\abs{\alpha_2}^2=1$; so
$\psi_k=\alpha_ke+\psi_k'$ with $\psi_k'\perp e$ has $\norm{\psi_1'}\norm{\psi_2'}=\abs{\alpha_1\alpha_2}$, and
$0=\abs h=\abs{a\bar\alpha_1\alpha_2+b\ip{\psi_1'}{\psi_2'}}\ge(a-b)\abs{\alpha_1\alpha_2}$.  Hence $\alpha_1\alpha_2=0$: one
of $\psi_1,\psi_2$ is orthogonal to $e$ and the other lies in $\C e$, and $\dim P\le2$, a contradiction.

So $\nu_c>0$.  We apply Definition~\ref{def:red-adapted}, Theorem~\ref{thm:red-square} and
Lemma~\ref{lem:red-KyFan} to the weight $D/\sqrt c$
(see \eqref{eq:red-normal}); the letters $D$, $a$, $b$, $c$ and $e$ keep their meaning.  The adapted bases are
unchanged, $\Upsilon_{D/\sqrt c}(\psi_1,\psi_2)=\norm X^2$, and from now on the numbers $h,r,s$, the matrix $B$ and
$R=r^2$, $S=s^2$, $H=\abs h^2$ are formed with $D/\sqrt c$.  Thus $G=1-\tau_\psi/c\ge0$, where $\tau_\psi$ is
formed with $D$, and $\nu=1-R-S=\nu_c/c>0$,
so Theorem~\ref{thm:red-square} applies.  Dividing the
proof of Lemma~\ref{lem:red-KyFan}(a) by $c$ gives, for unit vectors $v\in\C^2$,
\begin{equation}\label{eq:eq-Sp}
v^*S_pv=1-c^{-1}\Tr(\Pi_WD^2)\ge0,\qquad W=\operatorname{span}\{\psi_1,v_1e_2+v_2e_3\},
\end{equation}
(here $D$ is the original weight) and likewise for $S_q$ with $\psi_2$ and the $f_i$; so $B\in K^\C_2$.

We show $Rx_1=0$.  If $G=0$, then $x_1=(S_p)_{11}\ge0$ and $x_1=-\norm{B^*e_1}^2\le0$, so $x_1=0$.  Let $G>0$.
Then Theorem~\ref{thm:m2-D}~(iii) gives $F_\C(B)\succeq0$, so the terms on the right of \eqref{eq:red-square},
whose sum is $\norm{\mathsf Z}^2-\Upsilon_{D/\sqrt c}(\psi_1,\psi_2)=0$, are nonnegative and all vanish: $z_{33}=0$, the squares are $0$,
and $\mathbf x=(p;q)$ has $\mathbf x^*F_\C(B)\mathbf x=0$, hence $\mathbf x\in\ker F_\C(B)$.  Here $\mathbf x\ne0$, since for $p=q=0$ the vanishing
squares give $z_{32}=z_{23}=z_{22}=z_{11}=0$, hence $\mathsf Z=0$.  So $\det F_\C(B)=0$; both terms in
Theorem~\ref{thm:m2-D}~(i) are nonnegative (as in the proof of (ii) there), so $(G\Phi+Hd_p)(G\Phi+Hd_q)=0$, and
$\Phi\ge d_p,d_q\ge0$ gives $\Phi=0$.  As $\Phi=d_q+Rx_1$ with $d_q\ge0$ and $Rx_1\ge0$ on $K^\C_2$, we get $d_q=0$ and
$Rx_1=0$.  Now if $r=0$, then $D\psi_2\in\C\psi_1$ and $\psi_1\in\operatorname{span}\{\psi_2,e\}$.  If $r>0$, then $x_1=0$.
Since $x_1=(S_p)_{11}$, \eqref{eq:eq-Sp} with $v=e_1$ (the first standard basis vector of $\C^2$) gives
$\Tr(\Pi_WD^2)=c$ for $W=\operatorname{span}\{\psi_1,e_2\}=\operatorname{span}\{\psi_1,D\psi_2\}$, where $e_2$ is the
second vector of the adapted basis; by the remark before Lemma~\ref{lem:eq-2V}, $e\in W$, hence also $b\psi_2=D\psi_2-(a-b)\ip e{\psi_2}e\in W$.  In both cases $\dim P\le2$, a contradiction.
\end{proof}

\begin{proposition}[weights with $\lambda_2=\dots=\lambda_n$]\label{prop:eq-two-level}
Let $(X,Y)$ be an equality pair for $D=\diag(a,b,\dots,b)$, $a>b>0$, $n\ge3$.  Then $X,Y\in M(P)$ for
$P=\operatorname{span}\{\mathbf e_1,v\}$ with a unit vector $v\perp\mathbf e_1$.
\end{proposition}

\begin{proof}
By Lemma~\ref{lem:red-pure} and (2V),
\[
c\norm X^2\norm Y^2=\norm{T_D(X,Y)}^2\le\norm Y^2\max_\psi\Upsilon(\psi_1,\psi_2)\le c\norm X^2\norm Y^2,
\]
so $\Upsilon(\psi_1,\psi_2)=c\norm X^2$ for some unit vectors.  Apply Lemma~\ref{lem:eq-2V} to $X$.  In case (A),
Lemma~\ref{lem:eq-leak} gives $Y\in M(P)$.  Otherwise $X$ has the form (B); since $T_D(Y,X)=-T_D(X,Y)$, $(Y,X)$ is an
equality pair as well, and the same argument with $X$ and $Y$ exchanged gives a plane $P\ni e$ with $Y\in M(P)$, hence
$X\in M(P)$ by Lemma~\ref{lem:eq-leak}, unless $Y$ has the form (B) too.  In that case
$X=ew^*+w'e^*$ and $Y=eu^*+u'e^*$ with $w,w',u,u'\perp e$, and
\[
T_D(X,Y)=b\bigl(\ip w{u'}-\ip u{w'}\bigr)ee^*+a\,(w'u^*-u'w^*),
\]
whose two terms are orthogonal.  With $\norm X^2=\norm w^2+\norm{w'}^2$ and $\norm Y^2=\norm u^2+\norm{u'}^2$,
\begin{align*}
\abs{\ip w{u'}-\ip u{w'}}&\le\norm w\norm{u'}+\norm u\norm{w'}\le\norm X\norm Y,\\
\norm{w'u^*-u'w^*}&\le\norm{w'}\norm u+\norm{u'}\norm w\le\norm X\norm Y,
\end{align*}
by the Cauchy--Schwarz inequality, in $\C^n$ and then in $\R^2$, and $\norm{T_D(X,Y)}^2=c\norm X^2\norm Y^2$ forces
equality everywhere.  In $\R^2$ this gives $(\norm w,\norm{w'})=t(\norm{u'},\norm u)$ with $t>0$, so $w=0\iff u'=0$ and
$w'=0\iff u=0$; the first line gives $w\parallel u'$ and $u\parallel w'$, and the second (equality in the triangle
inequality for two rank-one matrices) $w'\parallel u'$ and $u\parallel w$, whenever the vectors involved are nonzero.  In
each of the three cases ($w,w',u,u'$ all nonzero; $w=u'=0$; $w'=u=0$) the nonzero ones among them are multiples of one
unit vector $v\perp e$, and $X,Y\in M(\operatorname{span}\{e,v\})$.
\end{proof}

For a general $D=\diag(t_1\ge t_2\ge\dots\ge t_n>0)$, $n\ge3$, $c=t_1^2+t_2^2$, let $J=\{k\ge3:t_k<t_2\}$; for
$\epsilon\in\{\pm1\}^J$ let $D_\epsilon$ be $D$ with the entries $t_k$, $k\in J$, replaced by $\epsilon_kt_2$, and
$D_+=\diag(t_1,t_2,\dots,t_2)$ the one with all $\epsilon_k=1$.

\begin{lemma}[transfer by convexity]\label{lem:eq-transfer}
If $(X,Y)$ is an equality pair for $D$, then $\norm{T_{D_\epsilon}(X,Y)}^2=c\norm X^2\norm Y^2$ and
$T_{D_\epsilon}(X,Y)=T_D(X,Y)$ for every $\epsilon$; in particular $(X,Y)$ is an equality pair for $D_+$.
\end{lemma}

\begin{proof}
Apply Lemma~\ref{lem:red-sign} to $s=(t_1,\dots,t_n)$ and $\beta=t_2$.  Its weights vanish whenever $\epsilon_k=-1$ for
some $k\ge2$ with $t_k=t_2$, and are positive otherwise, since $0<t_k\le t_2$; the remaining matrices are the
$D_\epsilon$, $\epsilon\in\{\pm1\}^J$, so $D=\sum_\epsilon w_\epsilon D_\epsilon$ with $w_\epsilon>0$ and $\sum_\epsilon w_\epsilon=1$.
The singular values of $D_\epsilon$ are $t_1,t_2,\dots,t_2$, so $\mathfrak n(D_\epsilon)^2\le c$ by
Theorem~\ref{thm:intro-norm}.  Hence
$\sqrt c\norm X\norm Y=\norm{\sum_\epsilon w_\epsilon T_{D_\epsilon}(X,Y)}\le\sum_\epsilon w_\epsilon\norm{T_{D_\epsilon}(X,Y)}\le\sqrt c\norm X\norm Y$.
So every $\norm{T_{D_\epsilon}(X,Y)}$ equals $\sqrt c\norm X\norm Y$, and by the strict convexity of the Frobenius norm all
$T_{D_\epsilon}(X,Y)$ are equal to their convex combination $T_D(X,Y)$.
\end{proof}

\begin{proof}[Proof of Theorem~\ref{thm:intro-eq}]
Sufficiency was shown after Proposition~\ref{prop:eq-n2}, and for $n=2$ necessity is that proposition.  Let
$n\ge3$.  A unitary change of basis preserves all conditions, so let $\omega$ be diagonal,
$D=\omega^{-1/2}=\diag(t_1\ge\dots\ge t_n)$, and let $(X,Y)=(A\omega^{1/2},B\omega^{1/2})$, an equality pair for $D$.  By
Lemma~\ref{lem:eq-transfer} it is an equality pair for $D_+$.  We first find a plane $P$ with $X,Y\in M(P)$.  If
$t_1>t_2$, Proposition~\ref{prop:eq-two-level} gives $P=\operatorname{span}\{\mathbf e_1,v\}$ with $v\perp\mathbf e_1$.  If
$t_1=t_2$, then $D_+=t_1I$ and $\norm{[X,Y]}^2=2\norm X^2\norm Y^2$; by the theorem of Cheng, Vong and Wenzel
\cite[Theorem~3.1]{CVW2010} (\cite[Theorem~2.7]{LuWenzel}) there is a unitary $U$
with $X=U(X_0\oplus0)U^*$, $Y=U(Y_0\oplus0)U^*$, $X_0,Y_0\in M_2$, and
$P=U(\C^2\oplus0)$ will do.

Let $\Pi=\Pi_P$ and $k\in J$.  Changing $\epsilon_k$ from $1$ to $-1$ changes $D_+$ by $2t_2E_{kk}$, so
$T_{E_{kk}}(X,Y)=0$ by Lemma~\ref{lem:eq-transfer}.  Suppose $u=\Pi\mathbf e_k\ne0$ and put $\hat u=u/\norm u$.  Since
$X,Y\in M(P)$, $0=T_{\Pi E_{kk}\Pi}(X,Y)=\norm u^2T_{\hat u\hat u^*}(X,Y)$.  Moreover
$\Pi D_+\Pi=\alpha\hat u\hat u^*+t_1u'u'^*$ for a unit vector $u'\in P$ orthogonal to $\hat u$ and some $\alpha\ge0$ (if
$t_1>t_2$, then $u=\ip v{\mathbf e_k}v$, $u'\in\C\mathbf e_1$ and $\alpha=t_2$; if $t_1=t_2$, then $\Pi D_+\Pi=t_1\Pi$).  Hence
$T_{D_+}(X,Y)=T_{\Pi D_+\Pi}(X,Y)=t_1T_{u'u'^*}(X,Y)$, and $\mathfrak n(u'u'^*)=\norm{u'u'^*}_{(2),2}=1$ (Theorem~\ref{thm:intro-norm}) gives
$\norm{T_{D_+}(X,Y)}^2\le t_1^2\norm X^2\norm Y^2<c\norm X^2\norm Y^2$, a contradiction.  Therefore $P\perp\mathbf e_k$ for
all $k\in J$, so $P\subseteq\operatorname{span}\{\mathbf e_k:k\notin J\}$, on which $D=D_+$.  Thus $P$ is invariant
under $D$ and $\omega$, and $\omega|_P$ has the eigenvalues $\lambda_1,\lambda_2$ ($D\mathbf e_1=t_1\mathbf e_1$ and $Dv=t_2v$
if $t_1>t_2$; $D|_P=t_1$ if $t_1=t_2$).  As $\Pi$ commutes with $\omega^{1/2}$, $A=\Pi A\Pi$ and $B=\Pi B\Pi$, and
Proposition~\ref{prop:eq-n2} for $\omega|_P$ gives $\varphi(A)=\varphi(B)=\varphi(A^*B)=0$.

Finally let $A,B$ be real and $\omega$ real symmetric, and let $P$ be any plane as in the theorem.  The sum $V$ of the
ranges of $A,A^*,B,B^*$ lies in $P$ and is spanned by real vectors.  If $\dim V=1$, then $V=\C u$ with a real unit vector
$u$ and $A=\alpha uu^*$ with $\alpha\ne0$, so $\varphi(A)=\alpha\,u^*\omega u\ne0$, a contradiction; hence $V=P$, and $P$ is
real (for $\omega=I$ compare \cite[footnote~24]{LuWenzel}).
\end{proof}

\begin{remark}[equality for general $Z$]\label{rem:eq-Z}
The argument also gives the cases of equality in \eqref{eq:intro-norm} when $s_2(Z)>0$.  Let $p,q\ge2$, let
$Z\in M_{q\times p}(\C)$ with $s_2(Z)>0$, and let $X,Y\in M_{p\times q}(\C)$ be nonzero.  Then equality holds in
\eqref{eq:intro-norm} if and only if there is a two-dimensional subspace $P\subseteq\C^p$, invariant under $Z^*Z$, on
which $Z^*Z$ has the eigenvalues $s_1(Z)^2,s_2(Z)^2$, such that, with $Q=ZP$, $X=\Pi_PX\Pi_Q$, $Y=\Pi_PY\Pi_Q$, and
the restrictions $X|_Q,Y|_Q\colon Q\to P$ and the inverse of $Z|_P\colon P\to Q$ are mutually orthogonal for the
trace inner product.  For $Z=\omega^{-1/2}$, $X=A\omega^{1/2}$ and $Y=B\omega^{1/2}$ this is Theorem~\ref{thm:intro-eq}.

To see this, write $Z=U\Sigma V$ as in the proof of Theorem~\ref{thm:intro-norm}.  Since
$T_\Sigma(VXU,VYU)=V(XZY-YZX)U$, and the padding used there preserves all the conditions (a plane on which
$Z_N^*Z_N$ has the positive eigenvalues $s_1(Z)^2,s_2(Z)^2$ lies in the range of $Z_N^*$, hence in $\C^p$), we may
assume $p=q=n$ and
$Z=\Sigma=\diag(t_1\ge\dots\ge t_n\ge0)$ with $t_2>0$.  In the proof of Theorem~\ref{thm:intro-eq}, the positivity of
$t_k$ for $k\ge3$ is used only in Lemma~\ref{lem:eq-transfer}, whose proof needs only $0\le t_k\le t_2$; the other
results are applied to the positive definite $D_+$.  So that proof gives a plane $P$, invariant under $\Sigma$, on which
$\Sigma$ has the eigenvalues $t_1,t_2$, with $X,Y\in M(P)$, and Proposition~\ref{prop:eq-n2} for the weight
$(\Sigma|_P)^{-2}$ on $P$ gives the orthogonality; the same proposition gives the converse.

For $s_2(Z)=0$ there are more equality pairs: for $Z=E_{11}\in M_3$, the pair $X=E_{11}$, $Y=E_{12}+E_{31}$ gives
$XZY-YZX=E_{12}-E_{31}$ and equality, although $Y^*$ has rank two and $\dim ZP\le1$ for every $P$.
\end{remark}

\section{Interpretation}\label{sec:interpretation}

Inequality \eqref{eq:14} is homogeneous of degree one in $\omega$ on both sides, so we may and do
assume $\Tr\omega=1$ in this section; throughout, $n\ge2$.  Then $\varphi(X)=\Tr(\omega X)$ is a faithful state on $M_n$,
and $\norm X_\omega^2=\varphi(X^*X)$.

\subsection{The GNS norm and inner derivations}
Let $L^2(M_n,\varphi)$ be $M_n$ with the inner product $\ip XY_\varphi=\varphi(X^*Y)$.  This is the
GNS space of $\varphi$: with $\xi_\varphi=\omega^{1/2}$ in the Hilbert--Schmidt space and
$\pi_\varphi(X)\xi=X\xi$, one has $\norm X_\omega=\norm{\pi_\varphi(X)\xi_\varphi}_2$.  For $A\in M_n$
let $\ad_A=[A,\cdot\,]$ be the inner derivation, viewed as an operator on $L^2(M_n,\varphi)$.
Theorem~\ref{cor:14} together with the example after Conjecture~\ref{conj:MKOC} says exactly that
\begin{equation}\label{eq:int-gns}
\sup_{A\ne0}\frac{\norm{\ad_A}_{L^2(\varphi)\to L^2(\varphi)}^2}{\norm A_\omega^2}
=\frac1{\lambda_1}+\frac1{\lambda_2},
\end{equation}
that is, the bilinear map $(A,B)\mapsto[A,B]$ on $L^2(M_n,\varphi)$ has norm
$(\lambda_1^{-1}+\lambda_2^{-1})^{1/2}$.  For the tracial state ($\omega=\frac1nI$) the right side is
$2n$, and \eqref{eq:int-gns} is \eqref{eq:BW} after the rescaling $\norm X_\omega^2=\frac1n\norm X_2^2$.

Under the unitary $X\mapsto X\omega^{1/2}$ of $L^2(M_n,\varphi)$ onto the Hilbert--Schmidt space,
$\ad_A$ becomes $L_A-R_{\sigma_{i/2}(A)}$, where $L_A\eta=A\eta$, $R_C\eta=\eta C$, and
$\sigma_{i/2}(A)=\omega^{-1/2}A\omega^{1/2}$ is the modular group $\sigma_t(A)=\omega^{it}A\omega^{-it}$ of
$\varphi$ at $t=i/2$; since $\norm\cdot_\omega$ is the left GNS norm, the right multiplication is twisted
by $\sigma_{i/2}$; this is why \eqref{eq:14}, although symmetric in $A$ and $B$, is not invariant under
$A\mapsto A^*$, $B\mapsto B^*$ (in general $\norm{A^*}_\omega\ne\norm A_\omega$).  The best
constant is not a function of the modular operator $\Delta_\varphi=L_\omega R_{\omega^{-1}}$ through its
norm: for tracial states $\Delta_\varphi=1$ for every $n$, while the constant is $2n$.  (By contrast, the
constant of \eqref{eq:15} is $1+\lambda_n/\lambda_1=1+\norm{\Delta_\varphi}$.)  It is a quantity of the
values of $\varphi$ on projections.  If $v$ is a nonzero partial isometry with $v^2=0$, then $p=vv^*$ and
$q=v^*v$ are orthogonal, $[v,v^*]=p-q$, and
\begin{equation}\label{eq:int-proj}
\frac{\norm{[v,v^*]}_\omega^2}{\norm v_\omega^2\,\norm{v^*}_\omega^2}=\frac{\varphi(p)+\varphi(q)}{\varphi(q)\,\varphi(p)}
=\frac1{\varphi(p)}+\frac1{\varphi(q)} .
\end{equation}
Since $\varphi(p),\varphi(q)\ge\lambda_1$ and $\varphi(p)+\varphi(q)=\varphi(p+q)\ge\lambda_1+\lambda_2$ ($p+q$ has rank at least $2$), the
right side is at most $\lambda_1^{-1}+\lambda_2^{-1}$ (the function $1/x+1/y$ is decreasing in each variable and convex
on the segment $x+y=\lambda_1+\lambda_2$, $\lambda_1\le x\le\lambda_2$), with equality for $v=e_1e_2^*$.
Thus Theorem~\ref{cor:14} says that the lower bound for the constant given by the pairs $v,v^*$ is sharp.

\begin{remark}[a KMS-twisted form]\label{rem:int-kms}
The adjoint of $\ad_A$ on $L^2(M_n,\varphi)$ involves the modular group at $t=-i$.  Indeed
$\ip{X}{YC}_\varphi=\Tr(\omega X^*YC)=\ip{X\omega C^*\omega^{-1}}{Y}_\varphi$, so the adjoint of $R_C$ is
$R_{\sigma_{-i}(C^*)}$ with $\sigma_{-i}(C)=\omega C\omega^{-1}$, and
$\ad_A^\dagger X=A^*X-X\sigma_{-i}(A^*)$.  Since $\norm{\ad_A^\dagger}=\norm{\ad_A}$, replacing $A$ by
$A^*$ shows that \eqref{eq:14} with the constant $\lambda_1^{-1}+\lambda_2^{-1}$ is equivalent to
\[
\norm{AX-X\sigma_{-i}(A)}_\omega^2\le\Bigl(\frac1{\lambda_1}+\frac1{\lambda_2}\Bigr)\varphi(AA^*)\,\varphi(X^*X)
\qquad(A,X\in M_n),
\]
with the same best constant.  The left side contains the twisted commutator of the KMS condition
$\varphi(AB)=\varphi(B\sigma_{-i}(A))$, and $\varphi(AA^*)$ is the square of the right GNS norm of $A$.
Equality holds for $A=v^*$ and $X=v$, where $v=e_1e_2^*$ is the extremal element of \eqref{eq:int-proj}.
\end{remark}

Since $[A,B]=[A-z-tB,\,B-w]$ for all scalars $z,w,t$, the inequality improves for free to a
Gram form, and hence to a variance form.  For $X\in M_n$ put $X_0=X-\varphi(X)$ and
$\Var_\varphi(X)=\norm{X_0}_\omega^2=\varphi(X^*X)-\abs{\varphi(X)}^2$.

\begin{corollary}\label{cor:int-var}
For every faithful state $\varphi=\Tr(\omega\,\cdot\,)$ on $M_n$ and all $A,B\in M_n$,
\begin{gather}
\norm{[A,B]}_\omega^2\le\Bigl(\frac1{\lambda_1}+\frac1{\lambda_2}\Bigr)
\Bigl(\Var_\varphi(A)\,\Var_\varphi(B)-\abs{\ip{A_0}{B_0}_\varphi}^2\Bigr),\label{eq:int-gram}\\
\norm{[A,B]}_\omega^2\le\Bigl(\frac1{\lambda_1}+\frac1{\lambda_2}\Bigr)\Var_\varphi(A)\,\Var_\varphi(B),
\label{eq:int-var}
\end{gather}
and the constant is best possible in both.
\end{corollary}

The bracket in \eqref{eq:int-gram} is the Gram determinant of $1,A,B$ in $L^2(M_n,\varphi)$, that is,
$\Var_\varphi(B)$ times the squared distance from $A$ to $\operatorname{span}\{1,B\}$.

\begin{proof}
If $\Var_\varphi(B)=0$, then $B_0=0$ since $\varphi$ is faithful, and both sides vanish.  Otherwise
apply \eqref{eq:14} to $A_0-tB_0$ and $B_0$, whose commutator is $[A,B]$, for $t\in\C$:
$\norm{[A,B]}_\omega^2\le(\lambda_1^{-1}+\lambda_2^{-1})\norm{A_0-tB_0}_\omega^2\Var_\varphi(B)$.  With
$g=\ip{B_0}{A_0}_\varphi$,
$\norm{A_0-tB_0}_\omega^2=\Var_\varphi(A)-\abs g^2/\Var_\varphi(B)+\Var_\varphi(B)\abs{t-g/\Var_\varphi(B)}^2$,
and $t=g/\Var_\varphi(B)$ gives \eqref{eq:int-gram}.  This is the argument of
\cite[Supplementary Information, Lemma~2]{KMO25}.  Inequality \eqref{eq:int-var} follows by dropping
$\abs{\ip{A_0}{B_0}_\varphi}^2$.  For the example $A=e_1e_2^*$, $B=A^*$ one has
$\varphi(A)=\ip{e_2}{\omega e_1}=0=\varphi(B)$ and $\ip{A_0}{B_0}_\varphi=\varphi(A^*A^*)=0$ since
$A^2=0$, so both \eqref{eq:int-gram} and \eqref{eq:int-var} are equalities.
\end{proof}

\begin{remark}[equality in the Gram and variance forms]\label{rem:eq-gram}
In Corollary~\ref{cor:int-var}, let $\Var_\varphi(B)>0$ and let the Gram determinant of $1,A,B$ be positive.  By its
proof, equality holds in \eqref{eq:int-gram} if and only if $(A_0-tB_0,B_0)$ with $t=\ip{B_0}{A_0}_\varphi/\Var_\varphi(B)$
gives equality in \eqref{eq:14}, and since $A_0-tB_0$ and $B_0$ are already orthogonal to $1$ and to each other,
Theorem~\ref{thm:intro-eq} says that this happens if and only if $A_0,B_0\in M(P)$ for an $\omega$-invariant plane $P$
on which $\omega$ has the eigenvalues $\lambda_1,\lambda_2$ (with $M(P)$ as in Section~\ref{subsec:m2-equality}).
Equality in \eqref{eq:int-var} holds if and only if, in addition, $\ip{A_0}{B_0}_\varphi=0$.
\end{remark}

When one of $A$, $B$ is normal, \eqref{eq:int-gram} is already in \cite{KMO25}: it is the observation
\cite[(16)]{KMO25} applied to the inequality \cite[(8)]{KMO25} for a normal matrix, proved in
\cite[Supplementary Information, Theorem~1]{KMO25}.  For Hermitian $A,B$ it is the uncertainty relation
\eqref{eq:int-KMO} below.

For $n=2$, \eqref{eq:int-gram} is an equality for all $A,B$: it is the identity of Proposition~\ref{prop:eq-n2}, which
for centred $A,B$ is \cite[(19)]{MKOC24}, written with the inner product $\ip\cdot\cdot_\varphi$.

\begin{remark}\label{rem:int-commutant}
Since $[A-Z,B_0]=[A,B]$ for every $Z$ in the commutant $\{B\}'=\{Z\in M_n:ZB=BZ\}$ of $B$, the proof of
\eqref{eq:int-gram} gives
$\norm{[A,B]}_\omega^2\le(\lambda_1^{-1}+\lambda_2^{-1})\operatorname{dist}_\varphi(A,\{B\}')^2\,\Var_\varphi(B)$,
where $\operatorname{dist}_\varphi$ is the distance in $L^2(M_n,\varphi)$, and the roles of $A$ and $B$ may be
exchanged.  As $\operatorname{span}\{1,B\}\subset\{B\}'$, this implies \eqref{eq:int-gram}; for $n=2$ the two
bounds coincide (the commutant of a non-scalar $2\times2$ matrix $B$ is $\operatorname{span}\{1,B\}$), and for
$n\ge3$ the new bound can be strictly smaller.  The constant is unchanged and is still attained by the
pair $A=e_1e_2^*$, $B=A^*$.  For comparison, the norm of $\ad_A$ as an operator on $B(\mathcal H)$ with
the operator norm is $2\operatorname{dist}(A,\C1)$ \cite{Stampfli1970}; the bound above is of the same
distance type in $L^2(M_n,\varphi)$, with the commutant of the other operator in place of the scalars.
\end{remark}

\begin{proposition}[the GNS form of Theorem~\ref{thm:intro-2V}]\label{prop:int-2V-gns}
For $A,B\in M_n$ put $J_A(B)=A^*B+B\,\sigma_{-i}(A^*)$, with $\sigma_{-i}(C)=\omega C\omega^{-1}$ as in
Remark~\ref{rem:int-kms}.  Then for all $A,B\in M_n$ with $B\ne0$,
\begin{equation}\label{eq:int-2V-gns}
\norm{[A,B]}_\omega^2+\operatorname{dist}_\varphi\bigl(J_A(B),\C B\bigr)^2
\le\Bigl(\frac1{\lambda_1}+\frac1{\lambda_2}\Bigr)\Var_\varphi(A)\,\norm B_\omega^2 ,
\end{equation}
and $\norm B_\omega^2$ cannot be replaced by $\Var_\varphi(B)$ \textup{(}for $B=1$ and $A=e_1e_2^*$ in an eigenbasis of
$\omega$, $J_A(1)=(1+\lambda_2/\lambda_1)\,e_2e_1^*$ is nonzero and orthogonal to $1$\textup{)}.
\end{proposition}

\begin{proof}
With $X=A\omega^{1/2}$, $Y=B\omega^{1/2}$ and $D=\omega^{-1/2}$ one has $\mathsf W=DX^*Y+YX^*D=J_A(B)\,\omega^{1/2}$.
Since $C\mapsto C\omega^{1/2}$ is a unitary of $L^2(M_n,\varphi)$ onto the Hilbert--Schmidt space,
$\norm{\mathsf W}_2=\norm{J_A(B)}_\omega$ and $\Tr(Y^*\mathsf W)=\ip B{J_A(B)}_\varphi$, so the last two terms on the
left of \eqref{eq:intro-2V} add up to $\operatorname{dist}_\varphi(J_A(B),\C B)^2$, and Theorem~\ref{thm:intro-2V} is
\eqref{eq:int-2V-gns} with $\norm A_\omega^2$ in place of $\Var_\varphi(A)$.  For a scalar $z$, $[A-z,B]=[A,B]$ and
$J_{A-z}(B)=J_A(B)-2\bar zB$, so the left side does not change when $A$ is replaced by $A_0=A-\varphi(A)$, and
$\norm{A_0}_\omega^2=\Var_\varphi(A)$.
\end{proof}

\begin{remark}\label{rem:int-2V-gns}
(a) Let $L_A$ and $R_A$ be the left and right multiplications by $A$ on $L^2(M_n,\varphi)$.  Then
$[A,B]=(L_A-R_A)B$, and $L_A^\dagger=L_{A^*}$ and $R_A^\dagger=R_{\sigma_{-i}(A^*)}$ (Remark~\ref{rem:int-kms}) give
$J_A(B)=(L_A+R_A)^\dagger B$.  Thus \eqref{eq:int-2V-gns} places the derivation $\ad_A=L_A-R_A$ next to the adjoint
of $L_A+R_A$, twice the Jordan multiplication by $A$; for a tracial state $J_A(B)=A^*B+BA^*$.  The twist
$\sigma_{-i}$ comes from the modular group of $\varphi$; the inner product is the GNS inner product $\varphi(X^*Y)$
throughout, not the KMS inner product $\Tr(\omega^{1/2}X^*\omega^{1/2}Y)$.

(b) For $\omega=I$ (Theorem~\ref{thm:intro-2V} does not need $\Tr\omega=1$) we have $D=I$, $c=2$ and
$\mathsf W=X^*Y+YX^*$, and Lemma~\ref{lem:red-Aud} is the identity
$\norm{[X,Y]}_2^2+\norm{\mathsf W}_2^2=\Tr\bigl((X^*X+XX^*)(Y^*Y+YY^*)\bigr)$ with which Audenaert begins his proof of
\eqref{eq:BW} \cite{Audenaert2010}.  Dividing by $4\norm Y_2^2$, \eqref{eq:intro-2V} becomes
\[
\Tr\Bigl(\rho\,\frac{X^*X+XX^*}2\Bigr)-\abs{\Tr(\rho X)}^2\le\frac12\norm X_2^2,\qquad
\rho=\frac{YY^*+Y^*Y}{2\norm Y_2^2}.
\]
Every density matrix is of this form ($Y=\rho^{1/2}$), and for all density matrices this is the variance bound
proved in \cite{Audenaert2010}.  So for scalar $\omega$, Theorem~\ref{thm:intro-2V} and
Proposition~\ref{prop:int-2V-gns} are equivalent to Audenaert's bound; what is new is the weight.
\end{remark}

\begin{remark}[a dual form]\label{rem:int-dual}
For $W,X,Y,Z\in M_n(\C)$, cyclicity of the trace gives
$\ip Z{XWY-YWX}=\Tr(YZ^*XW)-\Tr(XZ^*YW)=\ip{X^*ZY^*-Y^*ZX^*}W$.  Let
$\norm T^\circ=\max\{\abs{\ip ZT}:\norm Z_{(2),2}\le1\}$ be the dual norm of the Ky Fan $(2),2$ norm.  Taking the
maximum over $Z$ (resp.\ over $W$ with $\norm W_2=1$) in this identity shows that Theorem~\ref{thm:intro-norm} for
$p=q=n$ is equivalent to
\[
\norm{XWY-YWX}^\circ\le\norm W_2\norm X_2\norm Y_2\qquad(W,X,Y\in M_n(\C)).
\]
Since $\norm Z_{(2),2}\le\sqrt2\norm Z_{\mathrm{op}}$ and the trace norm is dual to the operator norm, the display implies
$\norm{XWY-YWX}_1\le\sqrt2\norm W_2\norm X_2\norm Y_2$, with equality for $n\ge2$, $W=E_{11}+E_{22}$, $X=E_{12}$,
$Y=E_{21}$.
\end{remark}

The constant involves $\lambda_1^{-1}$ itself, and the projection form \eqref{eq:int-proj} determines
which von Neumann algebras admit an analogue of \eqref{eq:14} with a finite constant.

\begin{proposition}\label{prop:int-vn}
Let $\mathcal M$ be a von Neumann algebra with a faithful normal state $\varphi$, put
$\norm X_\varphi^2=\varphi(X^*X)$, and let $c(\mathcal M,\varphi)\in[0,\infty]$ be the supremum of
$\norm{[A,B]}_\varphi^2/(\norm A_\varphi^2\norm B_\varphi^2)$ over nonzero $A,B\in\mathcal M$.
\begin{enumerate}
\item[\textup{(i)}] If $v\in\mathcal M$ is a nonzero partial isometry with $v^2=0$, then
$c(\mathcal M,\varphi)\ge\varphi(vv^*)^{-1}+\varphi(v^*v)^{-1}$.
\item[\textup{(ii)}] $c(\mathcal M,\varphi)<\infty$ if and only if
$\mathcal M\cong\mathcal A\oplus\bigoplus_{k=1}^NM_{n_k}(\C)$ with $\mathcal A$ abelian, $N<\infty$ and
$n_k\ge2$.
\item[\textup{(iii)}] In this case, let $\varphi=\Tr(\omega_k\,\cdot\,)$ on the summand $M_{n_k}(\C)$, with
$\omega_k$ positive definite and not normalised, and let $\lambda_1(\omega_k)\le\lambda_2(\omega_k)$ be its
two smallest eigenvalues.  Then
$c(\mathcal M,\varphi)=\max_k\bigl(\lambda_1(\omega_k)^{-1}+\lambda_2(\omega_k)^{-1}\bigr)$ \textup{(}$=0$ if
$N=0$\textup{)}, and, if $N\ge1$, equality holds in \textup{(i)} for some $v$.
\end{enumerate}
\end{proposition}

\begin{proof}
Part (i) is \eqref{eq:int-proj}, whose proof uses only $vv^*\perp v^*v$ and the faithfulness of $\varphi$.
Let $\mathcal M$ be as in (ii) and write $X=(X_{\mathcal A},X_1,\dots,X_N)$.  Then
$\norm X_\varphi^2=\norm{X_{\mathcal A}}_\varphi^2+\sum_k\norm{X_k}_{\omega_k}^2$ and $[A,B]=(0,[A_1,B_1],\dots,[A_N,B_N])$.
Theorem~\ref{cor:14} holds for $\omega$ not normalised, so with
$c_k=\lambda_1(\omega_k)^{-1}+\lambda_2(\omega_k)^{-1}$,
\[
\norm{[A,B]}_\varphi^2=\sum_k\norm{[A_k,B_k]}_{\omega_k}^2\le\sum_kc_k\norm{A_k}_{\omega_k}^2\norm{B_k}_{\omega_k}^2
\le\max_kc_k\,\norm A_\varphi^2\,\norm B_\varphi^2,
\]
and $v=e_1e_2^*$ in a summand with the largest $c_k$ gives equality in (i) with the value $\max_kc_k$.  This
proves (iii) and the implication $\Leftarrow$ in (ii).

For the standard facts on von Neumann algebras used in the rest of the proof (the partial isometry of the
polar decomposition of an element of $\mathcal M$ lies in $\mathcal M$; central projections and central supports)
see \cite[Proposition~II.3.14 and Section~V.1]{TakesakiI}.
Conversely, let $c(\mathcal M,\varphi)<\infty$, and let $z$ be the largest central projection with $\mathcal Mz$
abelian.  It suffices to show that $\mathcal M(1-z)$, which has no nonzero abelian direct summand, is
finite-dimensional.  If $f,g$ are orthogonal projections and $0\ne x\in f\mathcal Mg$, the partial isometry $v$ of the
polar decomposition of $x$ satisfies $vv^*\le f$ and $v^*v\le g$, hence $v^2=0$.  So every nonzero projection $f\le1-z$
dominates $vv^*$ for some $v\ne0$ with $v^2=0$: if $f\mathcal M(1-f)\ne0$, take $g=1-f$; otherwise $f$ is central,
$\mathcal Mf$ is not abelian and has a non-central projection $f'$, and $f'\mathcal M(f-f')\ne0$.  If $\mathcal M(1-z)$
were infinite-dimensional, it would contain mutually orthogonal nonzero projections $f_1,f_2,\dots$ (otherwise $1-z$ is
a finite sum of minimal projections $p_i$, each $p_i\mathcal Mp_j$ has dimension at most one, and $\mathcal M(1-z)$ is
finite-dimensional); then $\varphi(f_j)\to0$, and (i) gives $c(\mathcal M,\varphi)=\infty$.
\end{proof}

In particular, \eqref{eq:14} has no analogue with a finite constant for a faithful normal state on
$B(\mathcal H)$ with $\dim\mathcal H=\infty$, or on any factor of type II or III; on a II$_1$ factor with
its trace, matrix units of size $n$ give the ratio $2n$ in \eqref{eq:int-proj} for every $n$.

\begin{remark}[quaternionic coefficients]\label{rem:int-quaternion}
Over $\mathbb H$ the constant no longer sees $\lambda_2$.  Embed $M_n(\mathbb H)$ in $M_{2n}(\C)$ by
$X_1+X_2\mathbf j\mapsto\chi(X)=\bigl(\begin{smallmatrix}X_1&X_2\\-\overline{X_2}&\overline{X_1}\end{smallmatrix}\bigr)$, let
$\omega\in M_n(\mathbb H)$ be positive definite with quaternionic eigenvalues $\lambda_1\le\lambda_2\le\dots$, each of
multiplicity two in $\chi(\omega)$, and put $\norm X_{\omega}^2=\Tr(\chi(X)^*\chi(X)\chi(\omega))$.  Then the best
constant in $\norm{[A,B]}_\omega^2\le c\,\norm A_\omega^2\norm B_\omega^2$ is $2/\lambda_1$: Theorem~\ref{cor:14} for
$\chi(A),\chi(B),\chi(\omega)$ gives the upper bound, and in an eigenbasis $A=\mathbf iP_1$, $B=\mathbf jP_1$, with $P_1$ the
projection onto the first coordinate, attain it ($[A,B]=2\mathbf kP_1$).  This extremal pair spans a compact
$\mathfrak{su}(2)$ in the corner $P_1M_n(\mathbb H)P_1\cong\mathbb H$, whereas over $\R$ and $\C$ the extremal $v,v^*$ form
an $\mathfrak{sl}_2$-triple, which uses two eigenvalues.
\end{remark}

\subsection{Uncertainty relations}
Let $A,B$ be Hermitian (observables) and $C=[A,B]$; the state $\omega$ is denoted $\rho$ in
\cite{MKOC24,KMO25,SuzukiP12}.  Then $\Var_\varphi(A)=\Tr(\omega A^2)-(\Tr\omega A)^2$ is
the usual variance, and $\ip{A_0}{B_0}_\varphi=\operatorname{Cov}_\varphi(A,B)+\frac12\varphi(C)$, where
$\operatorname{Cov}_\varphi(A,B)=\frac12\varphi(A_0B_0+B_0A_0)$ is real and $\varphi(C)$ is purely
imaginary.  Hence \eqref{eq:int-gram} reads
\begin{equation}\label{eq:int-KMO}
\Var_\varphi(A)\Var_\varphi(B)-\operatorname{Cov}_\varphi(A,B)^2\ \ge\ \frac14\abs{\varphi(C)}^2
+\frac{\lambda_1\lambda_2}{\lambda_1+\lambda_2}\,\norm C_\omega^2 ,
\end{equation}
which is the uncertainty relation \cite[(4)]{KMO25}, and \eqref{eq:int-var} reads
\begin{equation}\label{eq:int-24}
\Var_\varphi(A)\,\Var_\varphi(B)\ \ge\ \frac{\lambda_1\lambda_2}{\lambda_1+\lambda_2}\,\norm C_\omega^2 ,
\end{equation}
the relation \cite[(24)]{MKOC24} derived there from the conjecture (the coefficient is denoted
$h$ in \cite{SuzukiP12}).  For Hermitian $A,B$ both relations are known (\cite{KMO25}, from the case of
\eqref{eq:14} with one normal matrix; \cite[Corollary~A and Remark~6.6]{SuzukiP12}, without \eqref{eq:14}); what
Theorem~\ref{cor:14} adds is \eqref{eq:int-gram} and \eqref{eq:int-var} for arbitrary matrices.

\begin{remark}[relation to \cite{SuzukiP12}]\label{rem:int-P12}
In \cite{SuzukiP12} the author studies commutator inequalities with a state-dependent weight from the
point of view of uncertainty relations.  For Hermitian $A,B$ it is shown there that
$\norm{[A,B]}_2^2\le\frac2{\lambda_1\lambda_2}\norm A_\omega^2\norm B_\omega^2$ \cite[Corollary~D(iii)]{SuzukiP12}; for
arbitrary complex $A,B$, Corollary~\ref{cor:16} gives the constant $\frac{\lambda_1+\lambda_2}{\lambda_1^2\lambda_2}$,
larger by the factor $\frac{\lambda_1+\lambda_2}{2\lambda_1}$ and attained by $A=\lambda_1E_{22}-\lambda_2E_{11}$,
$B=E_{21}$ \cite[\S2]{MKOC24}.  In both papers the best constant depends only on the two smallest eigenvalues of
the state, and the equality cases lie on the plane of the corresponding eigenvectors
(\cite[Theorem~D$'$]{SuzukiP12}, Theorem~\ref{thm:intro-eq}).  Both proofs reduce to extreme points, but of different
sets: the rank-two matrices $\pm i(xx^*-yy^*)$ in \cite{SuzukiP12}, where the bound follows from
${\wedge}^2\omega\ge\lambda_1\lambda_2$, and here the weights with $\lambda_2=\dots=\lambda_n$
(Proposition~\ref{prop:red-two-level}).  The two constants are paired: $\lambda_1\lambda_2$ is the smallest eigenvalue
of ${\wedge}^2\omega$, and $c=\lambda_1^{-1}+\lambda_2^{-1}$ is the largest eigenvalue of
$(\omega^{-1}\otimes I+I\otimes\omega^{-1})|_{{\wedge}^2\C^n}$ (both occur in \cite[Remark~2.8]{SuzukiP12}).
\end{remark}

\section{A concluding remark}\label{sec:conclusion}

Theorem~\ref{thm:intro-norm} bounds $\norm{XZY-YZX}_2$ by $\norm Z_{(2),2}\norm X_2\norm Y_2$, and the theorem of
Shi and Zhang \cite{ShiZhang} bounds it by $\sqrt2\norm Z_{\mathrm{op}}\norm X_{(2),2}\norm Y_2$; since
$\norm X_{(2),2}\le\norm X_2$ and $\norm Z_{(2),2}\le\sqrt2\norm Z_{\mathrm{op}}$, both would follow from the
following inequality.

\begin{conjecture}\label{conj:U}
For all $n\ge2$ and all $X,Y,Z\in M_n(\C)$,
\[
\norm{XZY-YZX}_2\le\norm Z_{(2),2}\,\norm X_{(2),2}\,\norm Y_2 .
\]
\end{conjecture}

The first of the two inequalities between norms is an equality when $X$ has rank at most two, in particular when
$n=2$, and the second when $s_1(Z)=s_2(Z)$, for instance when $Z=I_n$; in these cases Conjecture~\ref{conj:U} is
Theorem~\ref{thm:intro-norm} and the theorem of \cite{ShiZhang}, respectively, and so it holds.  It also holds when
$s_1(X)=s_2(X)$: if $X=U\diag(s)V$ with $U,V$ unitary and $s=(s_1(X),\dots,s_n(X))$, Lemma~\ref{lem:red-sign} with
$\beta=s_1$ writes $X$ as a convex combination of the matrices $s_1U\diag(1,1,\pm1,\dots,\pm1)V$, and the left side
is convex in $X$; for $X=\mu W$ with $W$ unitary, $XZY-YZX=\mu\,[WZ,YW^*]\,W$, and the refinement
$\norm{[A,B]}_2\le\sqrt2\norm A_{(2),2}\norm B_2$ of \eqref{eq:BW} due to Audenaert \cite{Audenaert2010} (see
\cite[(13)]{CVW2010}; it is also the case $Z=I_n$ of \cite{ShiZhang}), applied to $WZ$ and $YW^*$,
bounds its norm by $\sqrt2\,\abs\mu\,\norm Z_{(2),2}\norm Y_2$, which is
the right side because $\norm{\mu W}_{(2),2}=\sqrt2\,\abs\mu$.

\section*{Acknowledgements}
The author is grateful to the late Professor Takayuki Furuta, who introduced the author to operator
inequalities.  This research did not receive any specific grant from funding agencies in the public,
commercial, or not-for-profit sectors.

\section*{Data availability}
No data was used for the research described in the article.

\section*{Declaration of generative AI and AI-assisted technologies in the manuscript preparation process}
During the preparation of this work the author used Claude (Anthropic) in order to search for related
literature, to discuss possible approaches, to
write code for independent symbolic and numerical checks (Section~\ref{subsec:intro-computer}), and to improve the
language.  All references were checked against the original sources or their arXiv versions, and all proofs were verified and written out
by the author; none relies on computation.  After using this tool, the author reviewed and edited the content as
needed and takes full responsibility for the content of the published article.

\end{document}